\documentclass[12pt,a4paper]{article}
\usepackage{amsmath,amscd,amssymb,amsthm,amstext}

\input{xypic}

\theoremstyle{plain}
   \newtheorem{theorem}{Theorem}[section]
   \newtheorem{proposition}[theorem]{Proposition}
   \newtheorem{lemma}[theorem]{Lemma}
   \newtheorem{corollary}[theorem]{Corollary}
   \theoremstyle{definition}
   
   \newtheorem{definition}[theorem]{Definition}
   
   \theoremstyle{remark}
   
   \newtheorem{remark}[theorem]{Remark}

\newcommand{\FF}{{\mathbb F}}
\newcommand{\ZZ}{{\mathbb Z}}
\newcommand{\TT}{{\mathbb T}}

\newcommand{\K}{{\mathcal K}}
\newcommand{\F}{{\mathcal F}}
\newcommand{\Alg}{{{\mathcal A}lg}}
\newcommand{\R}{{\mathcal R}}
\newcommand{\sa}{{\mathcal A}}
\newcommand{\la}{{\Lambda}}
\newcommand{\Tor}{{\operatorname{Tor}}}
\newcommand{\tot}{{\operatorname{Tot}}}
\newcommand{\scat}{{\Delta}}
\newcommand{\sS}{{\mathcal S}}
\newcommand{\cS}{{c\mathcal S}}

\newcommand{\simp}{_\bullet}
\newcommand{\bisimp}{_{\bullet \bullet}}
\newcommand{\Hom}{{\operatorname{Hom}}}
\newcommand{\maps}{{\operatorname{map}}}
\newcommand{\mapc}{{\operatorname{Map}}}
\newcommand{\image}{{\operatorname{im}}}
\newcommand{\ex}{{E({\bf X})}}

\newcommand{\drf}{{\overline \Omega}}
\newcommand{\kerd}{{\mathcal Z}}
\newcommand{\imd}{{\mathcal B}}
\newcommand{\til}{{\widetilde \Omega}}
\newcommand{\lf}{{\mathcal L}}
\newcommand{\hd}{{\mathcal H}}
\newcommand{\Ot}{{\Omega_{tw}}}
\newcommand{\ea}{{\Lambda}}
\newcommand{\sdeg}{{s}}
\newcommand{\PS}{{P}}
\newcommand{\sig}{{\sigma}}
\newcommand{\bsig}{{\hat \sigma}}

\begin{document}
 
\title{A spectral sequence for string cohomology}
\author{Marcel B\" okstedt \& Iver Ottosen\thanks{The second author 
was supported by the University of Copenhagen and 
by the European Union TMR network ERB FMRX CT-97-0107: Algebraic 
$K$-theory, Linear Algebraic Groups and Related Structures.}}
\date{15 September 2004}

\maketitle

\begin{abstract}
Let $X$ be a 1-connected space with free loop space $\la X$. 
We introduce two spectral sequences converging
towards $H^*(\la X ;\ZZ /p)$ and $H^*((\la X)_{h\TT};\ZZ /p)$.
The $E_2$-terms are certain non Abelian derived
functors applied to $H^*(X;\ZZ /p)$. When  $H^*(X;\ZZ /p)$ 
is a polynomial algebra, 
the spectral sequences collapse for more or less trivial reasons.
If $X$ is a sphere it is a surprising fact that the spectral
sequences collapse for $p=2$.\\
\\
{\bf AMS subject classification (2000)}: 55N91, 55P35, 18G50
\end{abstract}

\section{Introduction}
Let $X$ be a space and let $\la X$ denote its 
free loop space. The circle group $\TT$ acts on $\la X$ by 
rotation of loops. The associated homotopy orbit space
$\la X_{h\TT}$  is sometimes called the string space.

For a manifold $X$, the free loop space has numerous
geometric applications. The most basic one is via the
Morse theory approach to the study of geodesic curves on $X$
\cite{GM}. But later the free loop space has also been used to
study diffeomorphisms. The main connection is through
Waldhausen's algebraic $K$-theory of spaces, 
the so called $A$-theory \cite{W}. One of the more refined 
versions of this connection relates pseudoisotopies
at a prime $p$ to the $p$-local space $TC(X,p)$.

One can define $TC(X,p)$ as the following homotopy  pullback:
$$
\begin{CD}
TC(X,p)^\wedge _p @>>> \Sigma^\infty (\Lambda X_+)^\wedge _p \\
@VVV @V1-\Delta_pVV \\
\Sigma^\infty (\Sigma (\Lambda X_{h\TT})_+)^\wedge _p @>Trf>>
\Sigma^\infty (\Lambda X_+)^\wedge _p
\end{CD}
$$
where $\Delta_p$ is the map which winds a loop $p$ times around itself
and $Trf$ is the $S^1$-transfer map. 
The space $TC(X,p)$ is an approximation to Waldhausen's
$A(X)$, which in turn gives a hold on the stable
pseudoisomorphism space of $X$. Even if the application
we have in mind is for differentiable manifolds,
it does not matter to $TC$ that $X$ is a manifold, 
as opposed to just a homotopy type. 

There is a third train of thoughts, inspired by
analogies with mathematical physics, especially 
with quantum field theory and string theory. In \cite{CS}
Sullivan and Chas introduced algebraic structures
relating  $H_*(\la X ;\ZZ )$ and $H_*(\la X_{h\TT} ;\ZZ )$. These
algebraic structures use that $X$ is a closed manifold. 
More precisely, the Thom class of the tangent bundle of 
$X$ plays an essential role.

The approach of this paper is homotopy theoretical.
We start with a homotopy type $X$, and try
to recover the modulo $p$ cohomology of $\la X$ and
of the Borel construction $\la X_{h\TT}$
by homotopy theoretical spectral sequences. 
We return to the construction of these spectral sequences
later in the introduction. The most essential properties are 
\begin{itemize}
\item It is derived from a cosimplicial space similar to the
cosimplicial space used to define the Adams spectral sequence 
for $X$.
\item The $E_2$ page of the spectral sequence is computed by 
non Abelian homological algebra in the sense of 
Andr\'e-Quillen homology.
\end{itemize}

One competing homotopy theoretical approach to the
cohomology of the Borel construction is  
the following three step method. Let us call it the
fibration method.
\begin{itemize}
\item Compute the cohomology of $\Omega X$ using
a spectral sequence (Serre or Eilenberg-Moore) 
belonging to the fibration $\Omega X \to PX \to X$.
\item Compute the cohomology of $\Lambda X$ using the
previous result and the fibration $\Omega X \to \Lambda X \to X$.
\item Compute the cohomology of the Borel construction
using the previous result and the fibration
$\la X \to \la X_{h\TT} \to B\TT $.
\end{itemize}

The Eilenberg Moore spectral sequence can be thought of
as the spectral sequence for the cohomology of a cosimplicial
space, just like our spectral sequence. But the particular
cosimplicial object is entirely different from ours. We use
a Postnikov decomposition of $X$, which is not visible
in the Eilenberg-Moore situation.

We feed two types of information into the machine, which are
not used by the fibration approach. First, we use that one
can explicitly compute the cohomology in the case when
$X$ is an Eilenberg-MacLane space. Secondly, we use non Abelian
homological algebra to keep track of how the pieces of our 
resolution fit together. In particular, a large part of the information 
about the $\TT$-action on $\Lambda X$ is internal to the
machine of non Abelian homological algebra.
So this part is taken care of already in the $E_2$ page.
 
This does not necessarily say that our method is better than the 
fibration method, however it does suggest that out method is
different. So there might be reasons to use both methods
simultaneously.

One possible drawback of our method is that non Abelian derived 
functors are hard to calculate. We have little use for a
spectral sequence with an incalculable $E_2$ page.

In order to show that this is not so bad, 
we do a few comparatively simple computations 
at the end of this paper. We show that at least for the
sphere $S^n$ we can solve the non Abelian homological algebra, 
and that the spectral sequence we get seems different from whatever 
comes out of the fibration method. We find it quite 
surprising and encouraging that in this case the 
spectral sequence converging to
$H^*(\la X_{h\TT} ;\FF_2 )$ collapses.
We do not actually prove this collapsing by methods internal to our 
spectral sequence. We quote the result which is known from  
other methods, and check by a counting argument that there 
is no room for differentials.

Because of the homotopy theoretical nature of our work,
it seems likely that we can use it to study $TC(X)$
for spaces $X$ with pleasant cohomology. We intend to
study this closer, but have not yet done so.

It seems more difficult, but potentially very 
profitable to compare our computations with
the Sullivan-Chas theory. We have not done this either yet.

What we have done, is that we have studied the spectral sequence
converging to $H^*(\Lambda X;\FF_2)$
in the very special case where $H^*(X;\FF_2)$ is a truncated polynomial
algebra on one generator. We can compute the relevant
$E_2$ page, and this makes it possible to
compute the Steenrod algebra action on 
$H^*(\Lambda X; \FF_2 )$ when $X$ is one of the projective 
spaces ${\mathbb C}P^n$, ${\mathbb H}P^n$ or the Cayley projective
plane $CaP^2$. The results led us to conjecture a stable splitting
of $\Lambda X$ for these spaces. We have later proved this 
splitting by unrelated methods.
 
The rest of the introduction is a more detailed description of 
our method. Consider the cohomology $H^*(X;\FF_p )$ as given.
The purpose of this paper is to study the cohomology of 
the free loop space and of its homotopy orbit space. 

In some cases, it is relatively easy to compute this cohomology.
For instance, suppose that $X$ is an Eilenberg-MacLane space. Then
there is a homotopy splitting $\la X \simeq X \times \Omega X$.
The space $\Omega X$ is also a Eilenberg-MacLane space, so that
the cohomology of $\la X$ is known.

The cohomology of the homotopy orbits $\la X_{h\TT}$ is 
more difficult to compute.
However, this is achieved in \cite{BO} and \cite{O2}.

The main idea of the present paper is to use these computations to study
the case of a general $X$. In essence, this application 
is done using a Postnikov decomposition of $X$. From our point of view, 
the simplest case is when $X$ is a product of Eilenberg-MacLane spaces, 
and correspondingly, the more $k$-invariants a space $X$ has, 
the more complicated it appears.
In particular, the spheres are very complicated spaces for this
approach.
 
Formally, we will study two spectral sequences converging towards 
the cohomology groups $H^*(\la X;\FF_p )$ and 
$H^*(\la X_{h\TT };\FF_p)$.
Both spectral sequences have origin in the Bousfield homology spectral
sequence \cite{B1}.

This is a remarkable spectral sequence that under fortunate circumstances
converges to the homology of the total space of a cosimplicial space. 

Let $X$ be a simply connected space. We re-write its Postnikov
tower as a cosimplicial space, whose total space is the
$p$-completion of $X$. This cosimplicial space is
the \emph{cosimplicial resolution} ${\bf R}X$ of $X$ with $R=\FF_p$.
Given this, we can form two cosimplicial spaces 
$\la {\bf R}X$ and $(\la {\bf R}X)_{h\TT}$ by applying
the functors $\la (-)$ and $\la (-)_{h\TT}$ in each codegree.
The total space of these new cosimplicial spaces are the completions of
$\la X$ respectively  $(\la X)_{h\TT}$. 
These cosimplicial spaces have associated Bousfield homology spectral
sequences $\{ \hat E^r \}$ and $\{ E^r \}$ respectively.

For 1-connected $X$ it is well known that $\{ \hat E^r \}$
converges strongly towards $H_*(\la X;\FF_p )$. We show that
$\{ E^r \}$ converges strongly towards $H_*(\la X_{h\TT};\FF_p )$
under the additional assumption that $H_*(X;\FF_p )$ is of finite type.

For the dual cohomology spectral sequences, 
$\{ \hat E_r \}$ and $\{ E_r \}$, 
we give an interpretation of the $E_2$ page.
The idea is that the $E_1$ page are given by the
cohomology of the respective functors (from spaces to spaces) 
applied  to the Eilenberg-MacLane spaces. This cohomology can,
according to \cite{BO},\cite{O1} and \cite{O2} be
written as certain functors $\drf$ respectively
$\ell$ (from algebras with a certain extra structure to algebras), 
applied to the cohomology of the Eilenberg-MacLane spaces.

This means that the $E_2$ page is the homology
of a chain complex, where the chains are given by these
functors applied to the cohomology of Eilenberg-MacLane spaces.
Since the cohomology of an Eilenberg-MacLane space
turns out to be a  free object, we can compute the $E_2$ pages
as derived functors.
 
To be precise, they are the non Abelian derived functor of $\drf$
applied to $H^*(X;\FF_p)$
respectively the non Abelian derived functor of $\ell$
applied to $H^*(X;\FF_p)$.
When $H^*(X;\FF_p)$ is a polynomial algebra the higher derived
functors vanish so 
the spectral sequences collapse at the $E_2$ page.

So far, the results are of a theoretical nature. As a
concrete example, we finally study the case $X=S^n$ and $p=2$. 
We develop homological algebra sufficient for computing the
relevant $E_2$ pages. 

For these spaces, there are other methods for computing 
$H^*(\la X ;\FF_p )$ and $H^*(\la X_{h\TT} ;\FF_p )$.

Comparing our $E_2$ pages with these results, we                          
show that for $X=S^n$ with $n\geq 2$ and $p=2$ 
the spectral sequences collapse at the $E_2$ pages.
 
We emphasize that this collapsing is not something to be expected 
a priori. Since spheres have complicated Postnikov systems, 
from the point of view of our spectral sequences, one would naively expect 
that these spectral sequence could have many nontrivial differentials.
So maybe the collapsing happens for a larger class of spaces?

Finally, we want to thank the referee of \cite{BO} for suggesting that we 
look at the Bousfield spectral sequence in this connection.

\section{Cosimplicial spaces with group actions}

In this section the category of simplicial sets is denoted 
$\sS$ and the category of cosimplicial spaces $\cS$. 
For $A,B\in \sS$ we let $\maps (A,B)= B^A$
denote the simplicial mapping space. We write ${\bf c} A$ for the constant
cosimplicial space with $({\bf c}A)^n=A$ for each $n$.

The category $\cS$ is a model category with weak equivalences, cofibrations and
fibrations as described in \cite{BK} X \S 4. The fibrations
are here defined in terms of matching spaces. By this definition 
it is clear that if $f:A \to B$ is a fibration in $\sS$ then
${\bf c} (f) : {\bf c}A \to {\bf c} B$ is a fibration in $\cS$.

The category $\cS$ is in fact a simplicial model category in the sense 
of \cite{Q} with ${\bf X}\otimes K \in \cS$, ${\bf X}^K \in \cS$ and 
$\mapc ({\bf X},{\bf Y})\in \sS$ defined as follows  
for $K\in \sS$ and ${\bf X},{\bf Y} \in \cS$:
\begin{eqnarray*} 
({\bf X}\otimes K)(\alpha ) &=& {\bf X}(\alpha ) \times K \\
({\bf X}^K)(\alpha ) &=& {\bf X}(\alpha )^K \\ 
\mapc ({\bf X},{\bf Y})_n &=& \Hom_{\cS} ({\bf X}\otimes \Delta^n ,{\bf Y})
\end{eqnarray*}
where $\alpha$ is a morphism in the simplicial category and 
$\Delta ^n= \Delta [n]\in \sS$ denotes the standard $n$-simplex.
In case $K$ is a simplicial group, this notation potentially clashes with
the usual notation for fixed points. In this paper, we are not going to
consider fixed points.

Let $\Delta$ be the cosimplicial space which in
codegree $n$ equals $\Delta^n$. We write $\Delta^{[m]}$ for the
simplicial $m$-skeleton and put $\Delta^{[\infty ]}=\Delta$.
By \cite{BK} X.4.3 we have that
$\Delta^{[m]}$ is a cofibrant cosimplicial space for each $0\leq m\leq \infty$.

The total space of a cosimplicial space ${\bf X}$ is defined as
$\tot {\bf X}=\mapc ({\bf \Delta}, {\bf X})$
If ${\bf X}$ is not fibrant, the total space might not
give you the ``right'' homotopy type. In this case, we
have to choose a fibrant replacement ${\bf Z}$ of ${\bf X}$, that is 
a weekly equivalent, fibrant cosimplicial space, and define 
$\overline \tot {\bf X}=\tot {\bf Z}$.

When the cosimplicial space has a group action one can choose an
equivariant fibrant replacement in the following sense:

\begin{lemma}
\label{lemma:equivariantreplacement}
Let $G$ be a simplicial group and ${\bf X}$ a cosimplicial $G$-space.
Assume that ${\bf X}^n$ is a fibrant simplicial set 
for each $n\geq 0$. Then there is a 
cosimplicial $G$-space $\ex$ such that both
$\ex$ and $\ex /G$ are fibrant cosimplicial spaces and such that 
the following diagram commutes: 
\begin{eqnarray}
\label{eqfibdiag}
\begin{CD}
EG \times {\bf X} @>\sim>> \ex \\ 
@VVV @VVV \\
EG \times_G {\bf X} @>\sim>>  \ex /G \\
\end{CD} 
\end{eqnarray}
Here the vertical maps are the obvious
quotient maps, and the horizontal maps are 
weak equivalences.
The map $\ex \to \ex /G$ is the pullback of the principal
$G$-fibration ${\bf c}EG \to {\bf c}BG$ over a fibration
$\ex /G \to {\bf c}BG$.
\end{lemma}

\begin{proof}
By the model category properties we can factor
the projection map  $EG\times_G {\bf X} \to {\bf c}BG$
as a composite $p\circ i$ where $i:EG\times_G {\bf X} \to {\bf Y}$ is
a cofibration which is simultaneously a weak equivalence, and
$p:{\bf Y} \to {\bf c}BG$ is  a fibration. $BG$ is a fibrant space by
\cite{GJ} Lemma I.3.5 so ${\bf c}BG$ is a fibrant cosimplicial space. Thus
${\bf Y}$ is fibrant.

We form the codegree wise pullback of $\pi : {\bf c}EG \to {\bf c}BG$ 
over $p$.
$$ 
\begin{CD}
\ex @= \ex @>{\bar p}>> {\bf c}EG \\
@VVV @V{\pi^p}VV @V{\pi}VV  \\
\ex /G @>\cong>> {\bf Y} @>p>> {\bf c}BG 
\end{CD} 
$$
The principal $G$-action (in the sense of \cite{May}) of $G$ on $EG$ gives
a principal $G$-action on $\ex^n$ for each $n$ and an isomorphism 
of cosimplicial spaces $\ex /G\cong {\bf Y}$ as written in the diagram.
By \cite{B1} Lemma 7.1 it follows that $\pi^p$ is a fibration so
$\ex$ is fibrant.

By the pullback property we can lift the map $i$ to a map 
$EG\times {\bf X} \to \ex$.
This constructs the missing map in the statement of the lemma.
In each codegree (\ref{eqfibdiag}) is a map of fibrations over $BG$ 
and we conclude that the lifting is also a weak equivalence. 
\end{proof}

\begin{theorem}
\label{toteq}
Let ${\bf X}$ be a fibrant cosimplicial space and $G$ a simplicial 
group. Then ${\bf X}^G$ is a cosimplicial $G$-space and we can
form its equivariant fibrant replacement $E({\bf X}^G)$. 
There is a natural map of fibrations of simplicial sets
for each $m$ with $0\leq m \leq \infty$:
$$\begin{CD}
(\tot_m {\bf X})^G @>>> EG\times_G (\tot_m {\bf X})^G @>>> BG \\
@V\sim VV @V\sim VV @V\cong VV \\
\tot_m (E({\bf X}^G)) @>>> \tot_m (E({\bf X}^G)/G) @>>> \tot_m ({\bf c}BG)
\end{CD}$$
The first and middle vertical maps are weak equivalences and the
right vertical map is an isomorphisms of simplicial sets.
\end{theorem}

\begin{proof}
Since ${\bf X}$ is fibrant each ${\bf X}^n$ is fibrant such that  
$({\bf X}^G)^n= ({\bf X}^n)^G$ is fibrant by \cite{May} Theorem 6.9.
Hence we can form $E({\bf X}^G)$.

By \cite{May} Definition 20.3 and Theorem 20.5 we have that the top 
vertical line in the diagram is a fiber bundle. 
By \cite{BK} X.5. SM7 and the fact that $\Delta^{[m]} \in \cS$ is cofibrant
we see that if $p:{\bf A} \to {\bf B}$ is a fibration in $\cS$ then
$\tot_m (p) : \tot_m {\bf A} \to \tot_m {\bf B}$ is a fibration in $\sS$.
In particular $\tot_m {\bf X}$ is fibrant since ${\bf X}$ is fibrant
and by \cite{May} Theorem 6.9 we have that $(\tot_m {\bf X})^G$ is fibrant.
Thus the top vertical line is a Kan fiber bundle and 
hence a fibration by \cite{May} Lemma 11.9.
The lower vertical line is $\tot_m$ of a fibration and hence a fibration.

There is a commutative diagram as follows:
$$ \begin{CD}
(\tot_m {\bf X})^G @>>> EG\times_G (\tot_m {\bf X})^G @>>> BG \\
@V\cong VV @Vf_m VV @V\cong VV \\
\tot_m ({\bf X}^G) @>>> \tot_m (EG\times_G {\bf X}^G) @>>> 
\tot_m ({\bf c}BG) \\
@V\sim VV @VVV @| \\
\tot_m (E({\bf X}^G)) @>>> \tot_m (E({\bf X}^G)/G) @>>> \tot_m ({\bf c}BG)
\end{CD} $$
The isomorphism $(\tot_m {\bf X})^G \cong \tot_m ({\bf X}^G)$ is one of 
the axiomatic isomorphisms in a simplicial model category. We examine it 
closer in order to define $f_m$.
A cosimplicial space is a diagram in $\sS$ and the axiomatic isomorphism comes 
from the corresponding isomorphism in the simplicial model category $\sS$.
For $A,B,C\in \sS$ this isomorphism is the composite
$$F: (A^B)^C \cong A^{B\times C} \cong A^{C\times B} \cong (A^C)^B$$

The following commutative diagram shows that $F$ is equivariant
with respect to actions of the monoid $C^C$.
$$\diagram
& C^C \times (A^B)^C \rrto^-{\circ} \dto & & (A^B)^C \dto \\
& C^C \times A^{B\times C} \rto^-{i_2\times 1} \dto 
&(B\times C)^{B\times C} \times A^{B\times C} \rto^-{\circ} 
& A^{B\times C} \dto\\
& C^C \times A^{C\times B} \rto^-{i_1\times 1} \dto 
&(C\times B)^{C\times B} \times A^{C\times B} \rto^-{\circ} 
& A^{C\times B} \dto \\
& C^C \times (A^C)^B \rto^-{i\times 1} 
& (C^C)^B \times (A^C)^B \rto^-{(\circ )^B} 
& (A^C)^B  
\enddiagram $$
For $Z\in \sS$ the action of $G$ on the mapping space $Z^G$ is 
defined by 
$$\begin{CD}
G\times Z^G @>ad(\mu )\times 1>> G^G \times Z^G @>\circ >> Z^G
\end{CD}$$
where $ad(\mu )$ denotes the adjoint of the product
$\mu :G\times G \to G$.
So taking $C=G$ in the above we see that $F$ is $G$-equivariant such 
that we have a map 
$$1\times_G F:EG\times_G (A^B)^G \to EG\times_G (A^G)^B$$
The composite
$$ \begin{CD}
EG\times (A^G)^B @>i \times 1>> (EG\times A^G)^B @>>> (EG\times _G A^G)^B 
\end{CD} $$
factors through $EG\times_G (A^G)^B$ and we compose with $1\times_GF$ to 
get a map
$$EG\times_G (A^B)^G \to (EG\times_G A^G)^B$$
The morphism $f_m$ in the theorem is codegree wise given
by this map. 

The lower part of the diagram is induced by (\ref{eqfibdiag}). 
The functor $(-)^K: \cS \to \cS$ where $K\in \sS$ preserves
fibrations as one sees from the right lifting property by taking
adjoints. Hence ${\bf X}^G$ is fibrant since ${\bf X}$ is fibrant.
By \cite{BK} X.5.2 we get a weak equivalence when applying $\tot_m$ to 
a weak equivalence between fibrant cosimplicial 
spaces. Thus the left vertical map is a weak equivalence.
The result follows.
\end{proof}

\section{Bousfield homology spectral sequences}

Let ${\bf X}$ be a fibrant cosimplicial space and let $A$ be an Abelian
group. In \cite{B1} Bousfield constructs a spectral sequence 
$\{ E^r({\bf X};A) \}$ with the homology of the  total space
$H_*(\tot {\bf X};A)$ as expected target. 

The precise convergence statement is as follows. Recall that there is
a tower of fibrations
$$\dots \to \tot_m {\bf X} \to \tot_{m-1} {\bf X} \to \dots
\to \tot_0 {\bf X}$$
with inverse limit $\tot {\bf X}$. Hence for each $n\geq 0$ there is a 
tower map
$$P_n ({\bf X}) : \{ H_n(\tot {\bf X};A) \}_{m\geq 0} 
\to \{H_n(\tot_m {\bf X};A) \}_{m\geq 0}$$
where the domain tower is constant. Let $A\otimes {\bf X}$ denote
the cosimplicial simplicial Abelian group with
$(A \otimes {\bf X})_t^m=A \otimes {\bf X}_t^m$ where 
$A\otimes S= \oplus_{x\in S}A$ for a set $S$.
Bousfield forms the double normalized complex and let $T(A\otimes {\bf X})$
denote its total complex. It is filtered by subcomplexes 
$F^mT(A\otimes {\bf X})$ and the quotient complex $T(A\otimes {\bf X}) /
F^{m+1}T(A\otimes {\bf X})$ is denoted $T_m(A\otimes {\bf X})$. 
A comparison map is defined
$$\Phi_n ({\bf X}) : \{ H_n(\tot_m {\bf X};A) \}_{m\geq 0} \to 
\{ H_nT_m(A\otimes {\bf X})\}_{m\geq 0}$$
and the following result is proved:

\begin{lemma}
\label{convergencecrit}
$\{ E^r({\bf X};A) \}$ converges strongly to
$H_*(\tot {\bf X};A)$ if and only if the tower map 
$\Phi_n({\bf X})\circ P_n({\bf X})$
is a pro-isomorphism for each $n$.
\end{lemma}

If $\Phi_n({\bf X})$ is a pro-isomorphism for each $n$ then ${\bf X}$ is
called an $A$-pro-convergent cosimplicial space and $\{ E^r({\bf X};A) \}$
is called pro-convergent.

We are interested in two special cases of this spectral sequence.
Let $R=\FF_p$ be the field on $p$ elements where $p$ is a fixed prime.
For a space $X$ we let ${\bf R}X$ denote the cosimplicial resolution of
$X$ in the sense of \cite{BK}. Note that $({\bf R}X)^n=R^{n+1}X$.
The free loop space on $X$ is by definition the simplicial mapping
space $\la X = \maps (\TT ,X)$ where we take $\TT =B\ZZ$.
By applying $\la$ codegree wise we get a cosimplicial space
$\la {\bf R}X$. We can also form the $\TT $ homotopy orbit space
codegree wise and get the cosimplicial space $(\la {\bf R}X)_{h\TT }$.
We are interested in the Bousfield homology spectral sequences for
these two spaces. 
As a corollary of \cite{B2} Proposition 9.7 we have

\begin{proposition}
\label{scla}
If $X$ is a 1-connected and fibrant space then 
$P_n(\la {\bf R}X)$ and $\Phi_n(\la {\bf R}X)$ are pro-isomorphisms for
each $n$ and the spectral sequence $\{ E^r(\la {\bf R}X; \FF_p) \}$
converges strongly to $H_*( \la (X_p^\wedge );\FF_p )\cong 
H_*(\la X;\FF_p )$.
\end{proposition}

\section{Strong convergence}
 
In this section we discuss convergence of the Bousfield homology
spectral sequence associated with $(\Lambda {\bf R } X)_{h\TT }$ where
$R=\FF_p$, the field on $p$ elements. We use $\FF_p$ coefficients 
everywhere unless stated otherwise. 

\begin{proposition}
\label{nilpotence}
If $X$ is a 1-connected space then $\la X$ and $\la X_{h\TT }$ are nilpotent
spaces. In fact we have $\pi_1(\la X)$- respectively 
$\pi_1(\la X_{h\TT })$-central series as follows for each $i\geq 1$:
\begin{align}
\pi_i (\la X) \supseteq \pi_i (\Omega X) \supseteq 0, & \label{cs1} \\
\pi_i (\la X_{h\TT }) \supseteq \pi_i (\la X) \supseteq 
\pi_i (\Omega X) \supseteq 0. & \label{cs2}
\end{align}
\end{proposition}

\begin{proof}
(\ref{cs1}) The fibration $\Omega X \to \la X \to X$ splits by the constant
loop inclusion $X\to \la X$. So
we have $\pi_i (\la X)\cong \pi_i  (\Omega X) \oplus \pi_i (X)$ for 
$i\geq 1$. Since the action of the fundamental group
is natural there is a commutative diagram

$$
\begin{CD}
\pi_1 (\Omega X) \times \pi_i (\Omega X) @>>> \pi_i (\Omega X) \\
@VVV @VVV \\
\pi_1 (\la X) \times \pi_i (\la X) @>>> \pi_i (\la X) \\
@VVV @VVV \\
\pi_1 (X) \times \pi_i (X) @>>> \pi_i (X)
\end{CD}
$$

We have $\pi_1 (\la X) \cong \pi_1 (\Omega X)$ since $X$ is simply 
connected. Further $\pi_1 (\Omega X)$ acts trivially on $\pi_i (\Omega X)$
since $\Omega X$ is an H-space. From the upper square we see that
the filtration (\ref{cs1}) is $\pi_1 (\la X)$-stable
and that the action on $\pi_i (\Omega X)$ is trivial. Since 
$\pi_1 (X)=0$ the lower square shows that the action on the 
quotient $\pi_i (\la X) / \pi_i (\Omega X)$ is trivial.

(\ref{cs2}) The fibration $\la X \to \la X_{h\TT } \to B\TT $ splits 
by a map constructed from a constant loop. 
So for $i\geq 1$ we have 
$\pi_i (\la X_{h\TT }) \cong \pi_i (\la X) \oplus \pi_i (B\TT  )$. 
Especially $\pi_1(\la X_{h\TT }) \cong \pi_1(\la X)$. 
By naturality there is commutative diagram

$$\begin{CD}
\pi_1 (\la X) \times \pi_i (\la X) @>>> \pi_i (\la X) \\
@VVV @VVV \\
\pi_1 (\la X_{h\TT }) \times \pi_i (\la X_{h\TT }) @>>> \pi_i (\la X_{h\TT }) \\
@VVV @VVV \\
\pi_1 (B\TT ) \times \pi_i (B\TT ) @>>> \pi_i (B\TT )
\end{CD} $$

From the upper square we see that the inclusion 
$\pi_i(\la X_{h\TT }) \supseteq \pi_i (\la X)$ 
is $\pi_1(\la X_{h\TT })$-stable. The lower square shows
that the action on the quotient $\pi_i(\la X_{h\TT }) / \pi_i(\la X)$ 
is trivial. The rest of the sequence (\ref{cs2}) has the desired properties
since (\ref{cs1}) is a $\pi_1(\la X)$-central series.
\end{proof}

\begin{proposition}
\label{proposition:procon}
If $X$ is a 1-connected space then the cosimplicial space 
$E(\la {\bf R} X)/\TT $ is $R$-pro-convergent. 
\end{proposition}

\begin{proof}
This is a consequence of \cite{B1} 3.3. 
Via the weak equivalences from Lemma \ref{lemma:equivariantreplacement}
we can use the filtrations from Proposition \ref{nilpotence} in each 
codegree. Then the quotients are $\pi_i ({\bf c}B\TT )$, $\pi_i ({\bf R}X)$ and
$\pi_{i+1} ({\bf R}X)$. Hence it suffices to show that when $n\leq 0$
the following holds for all $m\geq 0$:
\begin{equation}
\label{cohomotopy}
\pi^m \pi_{m+n} ({\bf c}B\TT )=0, \quad \pi^m \pi_{m+n} ({\bf R}X)=0, \quad 
\pi^m \pi_{m+n+1} ({\bf R}X)=0.
\end{equation}

Clearly $\pi^m \pi_{m+n} ({\bf c}B\TT )=0$ unless $m+n=2$ and 
$\pi^{2-n} \pi_2({\bf c}B\TT )=0$ since the differentials in the complex
$\pi_2 ({\bf c}B\TT )$ are alternating zeros and ones. 

By the proof of 6.1 in \cite{BK}, Ch. I and Proposition 6.3 in Ch. X
The following holds for any space $Y$: If $\widetilde H_i (Y; R)=0$ for
$i\leq k$ then $\pi^j \pi_i ({\bf R}Y)=0$ for $i\leq k+j$. So the last
two groups in (\ref{cohomotopy}) are also zero. 
\end{proof}

\begin{lemma}
\label{Rs}
Let $X$ be a 1-connected space with $H_*X$ of finite type.
Then $R_sX$ is 1-connected and $H_*R_sX$ is of finite type for each
$0\leq s < \infty$.
\end{lemma}

\begin{proof}
By \cite{BK} I.6.1 we have that $R_sX$ is 1-connected for each $s$.
Recall that $R(Y)$ is weakly equivalent to $\prod_{n=0}^\infty 
K(\tilde H_n (Y),n)$ for any space $Y$. So if $H_*Y$ is of finite
type then $H_*R(Y)$ is also of finite type and $\pi_iR(Y)= \tilde H_iY$
is finite for each $i$. Hence $\pi_i (({\bf R}X)^m)$ is finite for
each $i,m$. From \cite{S} Lemma 2.6 we see that $\pi_i(R_sX)$ is
finite for each $i,s$. By the Postnikov tower for $R_sX$ we conclude
that $H_*R_sX$ is of finite type for each $s$. 
\end{proof}

\begin{lemma}
\label{lim1}
Let $\dots \to C_*(2) \to C_*(1) \to C_*(0)$ be a sequence of maps 
of chain complexes. If for all $n$ and $m$ the group $C_n(m)$ is
finite, then there is an isomorphism  
$H_n(\lim C_*(m)) \cong \lim H_n(C_*(m))$
for all $n$.
\end{lemma}

\begin{proof}
This is a consequence of the $\lim^1$-sequence which can be 
found in e.g. \cite{Massey} Appendix A5.
\end{proof}

\begin{proposition}
\label{Piso}
Let $G$ be a simplicial group such that $H_n(BG)$ is finite
for all $n$. Let $\{ Z_m \}$ be a tower of $G$-spaces and put 
$Z_\infty =\lim Z_m$. Assume that $\{ H_*(Z_\infty ) \}_{m\geq 0} \to 
\{ H_*(Z_m) \}_{m\geq 0}$ is a pro-isomorphism and that $H_n(Z_m)$ 
is finite for all integers $n,m$. Then
$\{ H_*((Z_\infty)_{hG}) \}_{m\geq 0} \to 
\{ H_*((Z_m)_{hG}) \}_{m\geq 0}$ 
is also a pro-isomorphism.
\end{proposition}

\begin{proof}
We have Leray-Serre spectral sequences for $0\leq m \leq \infty$ as 
follows:
$$E^2(m)=H_*(BG;H_*(Z_m))\Rightarrow H_*((Z_m)_{hG}).$$
The tower map
$\{E^2_{i,j}(\infty ) \}_{m\geq 0} \to \{E^2_{i,j}(m) \}_{m\geq 0}$
is a pro-isomorphism for all $i$ and $j$ by the pro-isomorphism in 
the assumption, so 
$E^2_{i,j}(\infty)\cong \lim E^2_{i,j}(m)$. By the assumptions
on the homology of $BG$ and $Z_m$, the groups 
$E_{i,j}^2(m)$ with $m<\infty$ are all
finite so by Lemma \ref{lim1} we have 
$E^3_{i,j}(\infty) \cong \lim E^3_{i,j}(m)$.
By induction $E^r_{i,j}(\infty) = \lim E^r_{i,j}(m)$ for each $r$ and since
we have only finite filtrations 
$E^\infty_{i,j} (\infty) \cong \lim E^\infty_{i,j} (m)$.
Since $E^\infty_{i,j}(m)$ is finite for all $i,j,m$ it
follows that 
$\{E^\infty (\infty ) \}_{m\geq 0} \to \{E^\infty (m) \}_{m\geq 0}$
is a pro-isomorphism.
The result follows by the five lemma \cite{BK} III 2.7.
\end{proof}

\begin{theorem}
\label{strongcon}
If $X$ is a 1-connected fibrant space with $H_*(X;\FF_p)$ of finite type,
then the Bousfield spectral sequence
$\{ E^r(\la {\bf R}X_{h\TT };\FF_p) \}$ converges strongly to
$H_*(\la (X^\wedge_p)_{h\TT };\FF_p ) \cong 
H_*(\la X_{h\TT };\FF_p )$.
\end{theorem}

\begin{proof}
Let ${\bf Y}=\la {\bf R}X$.
The spectral sequence abuts to the homology of the total space
of a fibrant replacement of ${\bf Y}_{h\TT }$. 
We choose the fibrant replacement $E({\bf Y})/\TT $ 
from Lemma \ref{lemma:equivariantreplacement}. The total space of
this fibrant replacement is weakly equivalent to $\la (X_p^\wedge)_{h\TT }$
by Theorem \ref{toteq}. Thus the spectral sequence converges to 
the stated result if it converges. (A Leray-Serre spectral sequence 
argument shows that we can remove the $p$-completion inside the 
homology group.)

We have shown in Proposition \ref{proposition:procon} 
that the spectral sequence is pro-convergent.
Hence it suffices to show that $P_n(E({\bf Y})/\TT )$ or equivalently
$P_n({\bf Y}_{h\TT })$ is a pro-isomorphism.
By the Eilenberg-Moore spectral sequence and Lemma \ref{Rs} we see
that $H_*(\tot_s {\bf Y}) \cong H_*(\la R_s X)$ is of finite type for each
$0\leq s < \infty$. By Proposition \ref{Piso} and Proposition \ref{scla}
the result follows.
\end{proof}

We now change to cohomology. The dual of Proposition \ref{scla} and of
Theorem \ref{strongcon} is as follows:

\begin{theorem}
\label{dualspectralseq}
If $X$ is a 1-connected and fibrant space with $H_*X$ of finite
type then we have strongly convergent Bousfield cohomology spectral sequences 
\begin{align*}
& \hat E_r  \Rightarrow H^*(\la X), 
& & \hat E_2^{-m,t}=(\pi_m H^*(\la {\bf R}X))^t \\
& E_r \Rightarrow H^*((\la X)_{h\TT }), 
& & E_2^{-m,t} = (\pi_m H^*((\la {\bf R} X)_{h\TT }))^t .
\end{align*}
\end{theorem}

We are going to give a description of the $E_2$-terms as certain non
Abelian derived functors evaluated at $H^*X$. In the next section we
set up categories relevant for this purpose.

\section{The category $\F$ and the simplicial model category $s\F$}

For a fixed prime $p$ we let $\sa$ denote the mod $p$ Steenrod algebra and
$\K$ the category of unstable $\sa$-algebras. The category of
non-negatively graded unital $\FF_p$-algebras with the property that
$A^0$ is a $p$-Boolean algebra (ie. $x=x^p$ for all $x\in A^0$) is 
denoted $\Alg$. In \cite{O1}, \cite{O2} we defined a category
$\F$ with forgetful functors $\K \to \F \to \Alg$ as follows:

\begin{definition}
An object in $\F$ consists of an object $A$ in $\Alg$ which is 
equipped with an $\FF_p$-linear map 
$\lambda : A \to A$ with the following properties:
\begin{itemize}
\item $|\lambda x|=p(|x|-1)+1$ for all $x\in A$.
\item $\lambda x =x$ when $|x|=1$ and if $p$ is odd and $|x|$ is
even then $\lambda x=0$.
\item $\lambda (xy)=\lambda (x) y^p+x^p\lambda (y)$ for all $x,y \in A$. 
\end{itemize}
Furthermore $A$ is equipped with an $\FF_p$-linear map
$\beta : A \to A$ with the following properties:
\begin{itemize}
\item $|\beta x|=|x|+1$ for all $x\in A$.
\item $\beta \circ \beta =0$ and if $|x|=0$ then $\beta x=0$. 
\item $\beta (xy)=\beta (x)y + (-1)^{|x|}x\beta (y)$ for all $x,y \in A$. 
\end{itemize}
If $p=2$ we require that $\beta =0$.
A morphism $f:A\to A^\prime$ in $\F$ is an algebra homomorphism 
such that 
$f(\lambda x)=\lambda^\prime f(x)$ and $f(\beta x)=\beta^\prime f(x)$.
\end{definition}

\begin{remark}
For an object $K\in \K$ the map $\lambda :K \to K$ is defined by  
$\lambda x= Sq^{|x|-1}x$ when $p=2$ and $\lambda x= P^{(|x|-1)/2}x$ when
$p$ is odd and $|x|$ is odd. 
The map $\beta$ is the Bockstein operation when $p$ is odd.
\end{remark}

There is an obvious product on $\F$. There is also a coproduct. 
For two objects $A$ and $A^\prime$ in $\F$
the coproduct $A\otimes A^\prime$ is the tensor product of the underlying
objects in $\Alg$ equipped with maps $\lambda * \lambda^\prime$
and $\beta * \beta^\prime$ as follows; 
\begin{align*}
& \lambda * \lambda^\prime (x\otimes y)=
\lambda (x) \otimes y^p + x^p \otimes \lambda^\prime (y) \\
& \beta* \beta^\prime (x\otimes y)=
\beta (x) \otimes y + (-1)^{|x|} x\otimes \beta^\prime (y)
\end{align*}
In appendices in \cite{O1} and \cite{O2} we showed that $\F$ is complete
and cocomplete. It is well known that $\K$ and $\Alg$ also possess these 
properties. 

In the following $\R$ denotes any one of the categories $\K$, $\F$ or $\Alg$.
Let $n\FF_p$ denote the category of non-negatively graded 
$\FF_p$-vector spaces. The free functor $S_\R : n\FF_p \to \R$ 
is by definition the left adjoint of the forgetful functor $\R \to n\FF_p$. 
If $X$ is a non-negatively graded set we put
$S_\R (X)=S_\R (\FF_p \otimes X)$ where $\FF_p \otimes X$ is
the free graded $\FF_p$-vector space with basis $X$.
In particular we have free objects $S_\R (x_n)$ on one generator 
$x_n$ of degree $n$. 

\begin{remark}
Note that $S_\F (V)=S_\Alg (\bar V)$ where
\begin{align*}
\bar V & = V\oplus \bigoplus_{i\geq 1} \lambda^i V^{*\geq 2} & 
, \quad p=2 \\
\bar V & = V\oplus \beta V^{*\geq 1} \oplus 
\bigoplus_{i\geq 1, \nu \in \{ 0,1 \}} 
\beta^\nu \lambda^i (\beta V^{even,*\geq 2} \oplus V^{odd, *\geq 2}) &
, \quad p>2.
\end{align*}
\end{remark}

In the following we use \cite{Q} II.4 Theorem 4 to see that the 
category $s\R$ of simplicial objects in $\R$ 
is a simplicial model category. The arguments are standard but we
have included them anyhow.

We start by verifying that $\R$ has enough projectives. 
Recall that a morphism $f:X\to Y$ in a category $\mathcal D$ is called an
{\em effective epimorphism} if for any object $T$ and morphism 
$\alpha : X\to T$ there is a unique $\beta : Y \to T$ with 
$\beta \circ f= \alpha$ provided $\alpha$ satisfies the necessary condition 
that $\alpha \circ u = \alpha \circ v$ whenever $u,v : S \rightrightarrows X$
are maps such that $f\circ u = f\circ v$ 
(\cite{Q} II.4 proof of Proposition 2). 

\begin{proposition}
Let $f$ be an effective epimorphism in a category ${\mathcal D}$. Then
$f$ is an epimorphism. Furthermore if $f$ can be factored as $f=i\circ p$
where $i$ is a monomorphism then $i$ is an isomorphism.
\end{proposition}

\begin{proof}
Assume that $f$ is an effective epimorphism.  Let 
$r,s$ be two parallel arrows such that $r\circ f=s\circ f$.
Then for $\alpha = r\circ f$ we have $\beta \circ f = \alpha$ 
both for $\beta = r$ and $\beta =s$. So by uniqueness $r=s$. 
Thus $f$ is an epimorphism.

Assume that $f=i\circ p$ where $i$ is a monomorphism. 
If $f\circ u = f\circ v$ for two parallel arrows $u,v$ then 
$i\circ p \circ u = i \circ p \circ v$ and $p\circ u = p \circ v$ since
$i$ is a monomorphism. Hence there exists an arrow $j$ such that $p=j\circ f$.
Now $i \circ j \circ f= i\circ p =f$ which implies that $i\circ j = id$ 
since $f$ is an epimorphism.
Furthermore $i \circ j \circ i = id \circ i =i$ which implies that
$j\circ i = id$ since $i$ is a monomorphism.
\end{proof}

\begin{proposition}
A morphism in $\R$ is an effective epimorphism if and only if 
it is a surjection on underlying graded sets.
\end{proposition}

\begin{proof}
Any morphism $f:X \to Y$ may be factored as $X \to f(X) \to Y$ where 
the last map is clearly a monomorphism. So by the 
previous proposition we see that an effective epimorphism is 
surjective.

Assume that $f:X \to Y$ is a surjection and let $\beta : X \to T$ be
a map which satisfies $\beta \circ u = \beta \circ v$ whenever
$f\circ u = f\circ v$. For a given $x \in \ker f$ let $n=|x|$ and 
define $u,v : S_\R (x_n) \rightrightarrows  X$ by $u(x_n)=x$ and
$v(x_n)=0$. Then $x \in \ker \beta$ so we have $\ker f \subseteq \ker \beta$.
Now $\alpha : Y \to T$ with $\alpha (f(a))=\beta (a)$ is well defined 
and has $\alpha \circ f=\beta$.
\end{proof}

Recall that in \cite{Q} an object $P$ in a category $\mathcal D$ is called
{\em projective} if $\Hom_{D} (P,-)$ sends any effective epimorphism
to an surjection of hom-sets. 

\begin{proposition}
The following statements hold in the category $\R$:
\begin{enumerate}
\item $S_\R (V)$ is projective for any object $V$ in $n\FF_p$.
\item $\R$ has enough projectives.
\item $\{ S_\R (x_n) |n \geq 0 \}$ is a set of small projective 
generators.
\end{enumerate}
\end{proposition}

\begin{proof}
(1) By taking adjoints and applying the previous proposition we see that 
$S_\R (V)$ is projective. 
(2) Let $U:\R \to n\FF_p$ denote the forgetful functor and let $X$ be an
object in $\R$. The adjoint $\eta :S_\R (U(X)) \to X$ of 
$id_{U(X)}$ is surjective and hence an epimorphism. Thus $\R$ has enough 
projectives.
(3) The object $S_\R (x_n)$ is projective by (1). 
Since $\Hom_\R (S_\R (x_n),X)=X^n$ we have that $\Hom_\R (S_\R (x_n),-)$
commutes with filtered colimits so $S_\R (x_n)$ is small.
Finally, for two different morphisms $f,g :X\rightrightarrows Y$ 
there exist an $x\in X$ such that $f(x)\neq g(x)$. Hence the map
$S_\R (x_n) \to X$ with $x_n \mapsto x$ where $n=|x|$ 
separates $f$ and $g$ such that we have a set of generators as stated.
\end{proof}

We now turn to the category $s\R$ of simplicial objects in $\R$.
The homotopy groups of an object $R$ in $\R$ is
defined as the homology $\pi_*R=H_*(R,\partial )$ where $\partial$ is
the differential given by the alternating sums
$$\partial = \sum_{i=0}^n (-1)^id_i:
R_n \to R_{n-1}.$$ 
Especially $\pi_0(R)=R/(d_0-d_1)R$ and we have a morphism 
$\epsilon : R \to \pi_0(R)$ in $s\R$ given by projection where we
view $\pi_0(R)$ as a constant simplicial object. 

If $f:X\to Y$ is a morphism in $\R$ we can form the diagram
$$ \begin{CD}
X @>\epsilon>> \pi_0 X \\
@VfVV @V\pi_0 fVV \\
Y @>\epsilon>> \pi_0 Y
\end{CD} $$ 
One says that $f$ is {\em surjective on components} if the map 
from $X$ into the pullback 
$(f,\epsilon ): X \to Y \times_{\pi_0 Y} \pi_0 X$ 
is a surjection. Note that if $\pi_0(f)$ is an isomorphism then
$f$ is surjective on components if and only if $f$ is surjective.

\begin{proposition}
\label{smc}
There is a simplicial model category structure on $s\R$ as follows: 
\begin{itemize}
\item $f:X\to Y$ is a weak equivalence if 
$\pi_*f: \pi_*X \to \pi_*Y$ is an isomorphism.
\item $f:X\to Y$ is a fibration if it is surjective on components and
an acyclic fibration if it is both a fibration and a weak equivalence.
\item $f:X\to Y$ is a cofibration if for any commutative diagram
\begin{eqnarray}
\label{lift}
\begin{CD}
X @>>> A \\
@VfVV @VpVV \\
Y @>>> B
\end{CD} 
\end{eqnarray}
where $p$ is an acyclic fibration, there exist a map $Y \to A$ making 
both triangles commute.
\end{itemize}
The solution to the arrow diagram (\ref{lift}) is unique up to 
simplicial homotopy under $X$ and over $B$.
\end{proposition}

\begin{proof}
This is a special case of \cite{Q} II \S 4 Theorem 4. 
The uniqueness part follows from \cite{Q} II \S 2 Proposition 4.
\end{proof}

Note that the cofibrations are described in an indirect way. The concept 
of an almost free map make up for this weakness. See \cite{Q} II page 4.11
Remark 4 and the main source \cite{Miller1} \S 3, \cite{Miller2} \S 2 or
\cite{G}.

\begin{definition}
Let $\widetilde \Delta$ denote the subcategory of the simplicial category
$\Delta$ with objects $[n]=\{ 0,1,\dots ,n \}$ for $n\geq 0$ and 
morphisms the order preserving maps which sends $0$ to $0$.
An almost simplicial object in a category ${\mathcal C}$ is a functor
from $\widetilde \Delta^{op}$ to ${\mathcal C}$.
\end{definition}

\begin{definition}
A morphism $f:X\to Y$ in $s\R$ is called {\em almost free} if there
is an almost simplicial sub vector space $V$ of $Y$ such that
for each $n\geq 0$, the natural map $X_n\otimes S_\R (V_n) \to Y_n$
is an isomorphism.
\end{definition}

\begin{proposition}
(1) Almost free morphisms are cofibrations in $s\R$

(2) Any morphism $A \to B$ may be factored canonically and functorially as
$A\to X \to B$ where the first map is almost free and the second is an
acyclic fibration. 

(3) Any cofibration is a retract for an almost free map.
\end{proposition}

\begin{proof}
Similar to the one given in \cite{Miller1}. See also \cite{G}.
\end{proof}

\begin{definition}
A simplicial resolution of an object $A\in s\R$ is an acyclic 
fibration $P \to A$ in $s\R$ with $P$ cofibrant. 
An almost free resolution of $A$ is 
an acyclic fibration $Q \to A$ such that $\FF_p \to Q$ is almost free.
\end{definition}

Note that an almost free resolution is a resolution and that almost
free resolutions always exist by the above proposition.

\begin{theorem}
\label{cohcosres}
Let $X$ be a space with mod $p$ homology of finite type.
Let ${\bf R} X$ be the cosimplicial resolution of $X$. 
Then $H^*({\bf R}X)$ is an almost free resolution of $H^*X$ in 
each of the categories $\K$, $\F$ and $\Alg$.
\end{theorem}

\begin{proof}
This is a reformulation of well known results. We use \cite{GJ} VII 
Example 4.1 as a reference. Let $RX$ denote the simplicial 
$R$-module defined by applying the free $R$-module functor 
on each simplicial degree. Let $\eta : X \to RX$ be the map 
defined by $x\mapsto 1x$.

As in the \cite{GJ} one gets a cosimplicial space ${\bf R} X$ 
with $({\bf R}X)^n=R^{n+1}X$ and with augmentation 
$\eta : cX \to {\bf R} X$. Note that ${\bf R} X$ is a version of 
the Bousfield-Kan $R$-resolution of $X$ \cite{BK}. 
The homology $H_*(\eta ;R)=\pi_* (R\eta )$ is computed in \cite{GJ}
and taking the dual of this, we find that $\eta^*$ induces
an isomorphism as follows:
$$
\pi_s H^*({\bf R}X) \cong 
\begin{cases}
H^*X &, s=0 \\
0 &, s>0.
\end{cases}
$$
Thus $\eta^*$ is surjective and hence a fibration. Furthermore
$\eta^*$ is a weak equivalence. 

In order to show that $\FF_p \to H^*({\bf R} X)$ is almost free, 
we must find an almost simplicial sub vector space $V$ of 
$H^*({\bf R} X)$ such that $S_\R (V_n) \to H^*({\bf R} X)$ is an
isomorphism for each $n\geq 0$. 

As remarked in \cite{GJ} the cosimplicial maps $d^i$ for $i\geq 1$
and $s^i$ for $i\geq 0$ for ${\bf R}X$ are all morphisms of 
simplicial $R$-modules. Thus it suffices to show that 
$H^*({\bf R}X)^n$ is a free object in $\R$ for each $n\geq 0$.
But it is well known that $RX$ is homotopy equivalent to a
product of Eilenberg-MacLane spaces of the type $K(\FF_p , m)$.
The cohomology of such a product is a free object in $\K$ and 
hence also in $\F$ and $\Alg$.
\end{proof}

\section{Derived functors}

In this section $\R$ denotes any of the categories $\K$, $\F$ or $\Alg$.
We use the following notation for non Abelian derived functors:

\begin{definition}
The homology of an object $R$ in $s\R$ with coefficients
in a functor $E: \R \to \Alg$ is defined by
$$H_*(R;E) = \pi_* E(P)$$
where $P \to R$ is a simplicial resolution of $R$. By the uniqueness 
statement in Proposition \ref{smc} this homology theory is well defined 
and functorial in $R$.
\end{definition}

For an object $R \in \R$ we also write $R$ for the corresponding constant 
simplicial object in $s\R$. We are mainly interested in 
$H_*(R;E)$ when $R\in \R$. These homology groups have certain properties
which we now describe. 

Let $E$, $F$ and $G$ be functors from $\R$ to $\Alg$ 
with natural transformations $E \to F \to G$. Let $V: \Alg \to n\FF_p$ 
denote the forgetful functor to graded $\FF_p$-vector spaces.
If the sequence $0\to VE \to VF \to VG \to 0$ is short exact when
evaluated on any free object in $\R$ then we get a long exact sequence 
$$\dots \leftarrow H_i(R;E) \leftarrow H_i(R;F) \leftarrow H_i(R;G) \leftarrow
H_{i+1}(R;E) \leftarrow \dots .$$

The 0th homology group is sometimes given by the following result:
\begin{lemma}
\label{0th}
Define the category $\R^\prime $ as we defined the category $\R$
except that we do no longer require that objects are unital.
Let $F:\R^\prime \to \Alg^\prime$ be a functor. Assume that for every
surjective morphism $f:A\to B$ in $\R^\prime $ the following two conditions
hold:
\begin{enumerate}
\item $F(f):F(A) \to F(B)$ is surjective
\item $F(\ker f) \to F(A) \to F(B)$ is exact
\end{enumerate}
then $H_0(C;F)\cong F(C)$ for all objects $C$ in $\R$.
\end{lemma}

\begin{proof}
Let $P\to C$ be a simplicial resolution of $C$. From the normalized chain 
complex $N_*F(P)$ we see that $H_0(C;F)=F(P_0)/F(d_1)(\ker F(d_0))$. 

The maps $d_0,d_1: P_1 \to P_0$ are surjective by the simplicial identities. 
Let $i:\ker d_0 \to P_1$ denote the inclusion. By condition 2. we have that
$\ker F(d_0)=F(i)(F(\ker d_0))$. Thus
$$F(d_1)(\ker F(d_0)) = F(d_1)\circ F(i) (F(\ker d_0)).$$ 

There is a commutative diagram
$$\begin{CD}
\ker d_0 @>d_1^\prime >> d_1(\ker d_0) \\
@ViVV @VjVV \\
P_1 @>d_1>> P_0
\end{CD} $$
where $d_1^\prime$ denotes the restriction of $d_1$ and $j$ is the inclusion.
By this diagram $F(d_1)\circ F(i) = F(j)\circ F(d_1^\prime )$.
Furthermore $F(d_1^\prime )(F(\ker d_0))=F(d_1(\ker d_0))$ by condition 1.
So we have 
$$F(d_1)(\ker F(d_0))=F(j)\circ F(d_1^\prime )(F(\ker d_0))=
F(j)(F(d_1(\ker d_0)))$$
and 
$H_0(C;F)=F(P_0)/F(j)(F(d_1(\ker d_0))).$

Using condition 1. and 2. on the projection map $P_0\to P_0 /d_1(\ker d_0)$
we see that $H_0(C;F)\cong F(P_0/d_1(\ker d_0)) \cong F(C)$.
\end{proof}

The following result can sometimes be used to compute derived functors
of pushouts. We denote the pushout of a
diagram $A^\prime \leftarrow A \rightarrow A^{\prime \prime}$ in $s\R$ or $\R$ 
by $A^\prime \otimes_A A^{\prime \prime}$.

\begin{proposition}
\label{QuillenSS}
Let $E:\R \to \Alg$ be a functor. 

(1) If there is a natural isomorphism 
$E(A^\prime \otimes A^{\prime \prime})\cong 
E(A^\prime ) \otimes E(A^{\prime \prime})$ for objects 
$A^\prime ,A^{\prime \prime}$ in $\R$ 
then there is an isomorphism
$$H_*(B^\prime \otimes B^{\prime \prime};E) \cong 
H_*(B^\prime ;E) \otimes H_*(B^{\prime \prime };E) \text{ for }
B^\prime ,B^{\prime \prime } \in \R. $$

(2) Assume that there is a natural isomorphism 
$$E(A^\prime \otimes_A A^{\prime \prime})\cong 
E(A^\prime ) \otimes_{E(A)} E(A^{\prime \prime})$$
for diagrams 
$A^\prime \leftarrow A \rightarrow A^{\prime \prime}$ in $\R$.
Assume further that 
$B^\prime \leftarrow B \rightarrow B^{\prime \prime }$ is a diagram in $\R$ 
such that $\Tor_i^B(B^\prime , B^{\prime \prime })=0$ for $i>0$. 
Then there is a first quadrant spectral sequence as follows:
$$E^2_{i,j}=\Tor_i^{H_*(B;E)} (H_*(B^\prime ;E),H_*(B^{\prime \prime };E))_j 
\Rightarrow H_{i+j}(B^\prime \otimes_B B^{\prime \prime };E).$$
\end{proposition}

\begin{proof}
Let $P \to B$ be a simplicial resolution of $B$. By the factorization 
axiom we get a diagram 
$$
\diagram
& P^\prime \dto|>>\tip^\sim & P \lto|<<\tip \rto|<<\tip \dto|>>\tip^\sim 
& P^{\prime \prime } \dto|>>\tip^\sim \\
& B^\prime & B \lto \rto & B^{\prime \prime } 
\enddiagram 
$$
where the vertical maps are acyclic fibrations and the upper horizontal
maps are cofibrations as indicated. 
Since $\FF_p \to P$ is a cofibration and cofibrations are stable under
composition we see that $P^\prime \to B^\prime$ and
$P^{\prime \prime } \to B^{\prime \prime }$ are simplicial resolutions.

Now form the map of pushouts 
$f: P^\prime \otimes_P P^{\prime \prime } \to 
B^\prime \otimes_B B^{\prime \prime }$ and consider the corresponding
map of derived tensor products in the sense of \cite{Q} II \S 6:
$$Lf :  P^\prime \otimes_P P^{\prime \prime } 
= P^\prime \otimes_P^L P^{\prime \prime } \to   
B^\prime \otimes_B^L B^{\prime \prime }.$$
By \cite{Q} II \S 6 Theorem 6 there are second quadrant spectral sequences
\begin{align*}
\Tor_i^{\pi_*P} (\pi_*P^\prime , \pi_*P^{\prime \prime})_j 
\Rightarrow \pi_{i+j} (P^\prime \otimes_P P^{\prime \prime}) \\
\Tor_i^{\pi_*B} (\pi_*B^\prime , \pi_*B^{\prime \prime})_j 
\Rightarrow \pi_{i+j} (B^\prime \otimes_B^L B^{\prime \prime})
\end{align*}
The above diagram gives a map of spectral sequences which 
is an isomorphism at the $E_2$-terms. Hence $Lf$ is a weak equivalence.
By the Corollary following Quillen's Theorem 6 we have that
$B^\prime \otimes_B^L B^{\prime \prime} \to 
B^\prime \otimes_B B^{\prime \prime}$ is a weak equivalence. Thus
$f$ is itself a weak equivalence.

Since $f$ is surjective it is a fibration. Since the pushout of a cofibration
is a cofibration $P^\prime \to P^\prime \otimes_P P^{\prime \prime }$
is a cofibration and thus the domain of $f$ is cofibrant. So $f$ is
a simplicial resolution.

For the proof of (1) take $B=\FF_p$ and apply $E$ codegree wise.
The result follows by the Eilenberg-Zilber theorem.
For the proof of (2) apply $E$ codegree wise. The result follows by 
Quillen's Theorem 6.
\end{proof}

If one knows that the higher derived functors vanish on a certain
class of objects, they can be used to compute derived functors
by the following result.

\begin{theorem}
\label{genres}
Let $E:\R \to \Alg$ be a functor and let $A\in \R$.
Assume that $Q \tilde \twoheadrightarrow A$ is an acyclic fibration in 
$s\R$ and that 
$$
H_i(Q_j;E)= \begin{cases}
E(Q_j) & , i=0 \\
0 & , i>0.
\end{cases} 
$$
Then $H_*(A;E)\cong \pi_* E(Q)$.
\end{theorem}

\begin{proof}
We have shown that $s\R$ is a simplicial model category. So $ss\R$ is
a simplicial model category by the Reedy structure \cite{GJ} VII 2.13.
A fibration in $ss\R$ is especially a level fibration and a cofibration
is especially a level cofibration by \cite{GJ} VII 2.6. A weak equivalence
is a level weak equivalence by definition. 

We use a dot to denote a simplicial direction in the following.
Let $cQ\bisimp$ denote the object in $ss\R$ defined by $(cQ)_{ij}=Q_j$ for 
all $i$. Let $P\bisimp$ be a resolution of $cQ\bisimp$ ie.
$(\FF_p)\bisimp \rightarrowtail P\bisimp 
\tilde \twoheadrightarrow cQ\bisimp$. 

We have that 
$(\FF_p)\simp \rightarrowtail P_{i \bullet} \tilde \twoheadrightarrow Q\simp$ 
for each $i$ by the above. By composition with the acyclic fibration
$Q\simp \tilde \twoheadrightarrow A$ we see that 
$P_{i \bullet}$ is a resolution of $A$.

So the horizontal homotopy of $E(P\bisimp )$ is given by  
$\pi_j^h E(P_{i \bullet}) = H_j(A;E)$. We apply vertical homotopy on
this and obtain 
$$
\pi_i^v \pi_j^h E(P\bisimp ) \cong 
\begin{cases}
H_j(A;E) & , i=0 \\
0 & , i>0.
\end{cases}
$$ 

We also have that $P_{\bullet j}$ is a resolution of $Q_j$ 
for each $j$. So $\pi_i^v E(P_{\bullet j}) \cong H_i(Q_j;E)$ which equals 
$E(Q_j)$ for $i=0$ and equals $0$ for $i>0$.
We apply horizontal homotopy on this and obtain
$$\pi_j^h \pi_i^v E(P\bisimp ) \cong 
\begin{cases}
\pi_j E(Q\simp ) & , i=0 \\
0 & , i>0.
\end{cases}
$$
Thus both spectral sequences associated with $E(P\bisimp )$ collapse
and the result follows.
\end{proof}

\section{The $E_2$ terms seen as derived functors}

In \cite{BO}, \cite{O1} and \cite{O2} we introduced a functor 
$\drf :\F \to \Alg$ as follows:

\begin{definition}
$\drf (R)$ is the quotient of the free graded commutative and unital 
$R$-algebra on generators $\{ dx | x\in R \}$ of degree $|dx|=|x|-1$, 
modulo the ideal generated by the elements
\begin{align*}
& d(x+y)-dx-dy, & & d(xy)-d(x)y-(-1)^{|x|}xd(y), \\
& d(\lambda x)-(dx)^p, & & d(\beta \lambda x).
\end{align*}
There is a differential $d:\drf (R) \to \drf(R)$ given by
$d(x)=dx$ for $x\in R$.
\end{definition}

Note that for $p=2$ the Bockstein is trivial so here the functor $\drf$ 
is the same as the functor which we originally denoted $\Omega_\lambda$. 

It was shown that there is a lift to a 
functor $\drf : \K \to \K$ and that this lift is nothing but
Lannes' division functor $(-:H^*(\TT ))_\K$. In particular 
there is a morphism $\drf (H^*X) \to H^*(\Lambda X)$ for any space $X$
which is an isomorphism when $H^*X$ is a free object in $\K$.

An other functor $\ell: \F \to \Alg$ was also introduced in 
\cite{BO}, \cite{O1}, \cite{O2} as follows: 

\begin{definition} Let $p=2$ and let $A$ be an object in $\F$.
The $\FF_2$-algebra $\ell (A)$ is the quotient of the free 
graded commutative $\FF_2$-algebra on generators
$$\phi (x), q(y), \delta (z) , u \text{ for } x,y,z \in A$$
of degrees $|\phi (x)|=2|x|$, $|q(x)|=2|x|-1$, $|\delta (x)| =|x|-1$
and $|u|=2$, by the ideal generated by the elements

\begin{align*}
& \phi (a+b)+\phi (a) +\phi (b) , \\
& \delta (a+b)+\delta (a)+ \delta (b), \\
& q(a+b)+q(a)+q(b)+\delta (ab) \\
& \delta (xy)\delta (z)+\delta (yz)\delta (x)+\delta (zx)\delta(y) ,\\
& \phi (xy)+ \phi (x) \phi (y) +uq(x)q(y), \\
& q(xy)+q(x)\phi (y) +\phi (x)q(y), \\
& \delta (x)^2 +\delta (\lambda x), \\
& q(x)^2+\phi (\lambda x)+ \delta (x^2\lambda x), \\
& \delta (x)\phi (y) +\delta (xy^2), \\
& \delta (x) q(y) +\delta (x\lambda y)+\delta (xy)\delta (y), \\
& u\delta (x), 
\end{align*}
where $a,b,x,y,z$ are homogeneous elements in $A$ with $|a|=|b|$.
\end{definition}

\begin{definition}
Let $p$ be an odd prime and let $A$ be an object in $\F$. 
The $\FF_p$-algebra $\ell (A)$ is the quotient of the  
free graded commutative $\FF_p$-algebra on generators 
$$\phi (x), q(y), \delta (z), u \text{ for } x,y,z \in A$$ 
of degrees $|\phi (x)|=p|x|-\sig (x)(p-1)$, 
$|q(y)| = p|y|-1-\sig(y)(p-3)$, $|\delta (z)|= |z|-1$ and $|u|=2$
by the ideal generated by the elements 
\begin{align*}
& \phi (a+b)-\phi (a)-\phi (b)+\sig (a)
\sum_{i=0}^{p-2}(-1)^i \delta (a)^i \delta (b)^{p-2-i}\delta (ab), \\
& \delta (a+b)-\delta (a) -\delta (b), \\
& q(a+b)-q(a)-q(b)+\bsig (a)
\sum_{i=1}^{p-1} (-1)^i \frac 1 i \delta (a^ib^{p-i}), \\ 
& (-1)^{\sig (x) \bsig(z)} \delta (x)\delta (yz) + 
(-1)^{\sig (y) \bsig (x)} \delta (y)\delta (zx) +
(-1)^{\sig (z) \bsig (y)} \delta (z)\delta (xy), \\
& \phi (xy)-(-u^{p-1})^{\sig (x) \sig (y)} \phi (x) \phi (y), \\
& q(xy)-(-u^{p-1})^{\sig (x) \sig (y)}
(u^{\sig (y)}q(x)\phi (y)+(-u)^{\sig (x)}\phi (x) q(y)), \\
& q(x)^p - u^{p-1} q(\lambda x) -\phi (\beta \lambda x), \\
& \delta (x)\phi (y)-\delta (xy^p)
-\delta (x\lambda y)+\delta (xy)\delta (y)^{p-1}, \\
& \delta (x) q(y)-\delta (xy^{p-1})\delta (y)-
\delta (x\beta \lambda y), \\
& \delta (x)u, \\
& q(\beta \lambda x), \\
& \delta (x^p), 
\end{align*}
where $a, b, x, y, z \in A$ with $|a|=|b|$. Furthermore, 
$\sig (x)=1$ for $|x|$ odd, $\sig (x)=0$ for $|x|$ even and
$\bsig (x) =1-\sig (x)$.
\end{definition}

The functor $\ell$ also lifts to an endofunctor on $\K$ and it comes with 
a natural morphism $\ell (H^*X) \to H^*(\Lambda X_{h\TT })$ which is
an isomorphism if $H^*X$ is a free object in $\K$. For details on
this see \cite{BO}, \cite{O1}, \cite{O2}.

Via Theorem \ref{cohcosres} we can now restate Theorem 
\ref{dualspectralseq} in an appropriate form. 

\begin{theorem}
\label{mainss}
If $X$ is a 1-connected and fibrant space with $H_*X$ of finite
type then we have strongly convergent Bousfield cohomology spectral sequences 
$\hat E_r  \Rightarrow H^*(\la X)$ and $E_r \Rightarrow H^*(\la X_{h\TT })$
with the following $E_2$ terms:
$$
\hat E_2^{-m,t} \cong H_m(H^*(X);\drf)^t \text{ and }
E_2^{-m,t} \cong H_m(H^*(X);\ell )^t.
$$
\end{theorem}

We now introduce other functors in order to study the derived functors
of $\ell$. Recall that the functors $\lf$ and $\til$ from $\F$ to $\Alg$ 
are defined by $\lf (R) = \ell (R)/(u)$ and 
$\til (R) = \lf (R)/ (\delta (x)|x\in R)$.

\begin{proposition}
\label{0thapl}
For each object $R\in \F$ there are isomorphisms as follows:
$H_0(R;\drf)\cong \drf (R)$, $H_0(R;\til )\cong \til (R)$ and 
$H_0(R;\lf )\cong \lf (R)$.
\end{proposition}

\begin{proof}
We use Lemma \ref{0th} to prove this.
By their definitions we may consider $\drf$, $\til$ and $\lf$ as functors from 
$\F^\prime$ to $\Alg^\prime$.

Let $A$ be an object in $\F^\prime$ and let $I\subseteq A$ be an ideal.
We must verify condition 1. and 2. in Lemma \ref{0th} for these functors
where $f: A\to A/I$ is the natural projection. We do this for the functor
$\lf$. The verification for the other functors is similar 
but easier.

The map $\lf (f)$ is surjective with kernel 
$$J=(\phi (x)-\phi (y),q(x)-q(y),\delta (x)-\delta (y)|x-y\in I) \subseteq 
\lf (A) $$ 
so $\lf (A)/J \cong \lf (A/I)$. We must check that $\lf (I)= J$. 

The inclusion $\lf (I) \subseteq J$ holds
since $\phi (0)=q(0)=\delta (0)=0$.

For the inclusion $\lf (I) \supseteq J$ assume first that
$p=2$. Since $\phi$ and $q$ are additive we have that
$\delta (x)-\delta (y)$ and $\phi (x)-\phi (y)$ lie in $\lf (I)$.
Further $q(x)-q(y)=q(x-y)+\delta (xy)$ but $\delta (xy)=\delta (x(y-x))$
so also $q(x)-q(y)\in \lf (I)$. Thus the inclusion holds.

For $p$ odd $\delta$ is additive, $\phi$ is additive on elements
of even degree and $q$ is additive on elements of odd degree. 
For $|x|=|y|$ odd we have
$$\phi (x)-\phi (y)=\phi (x-y)+
\sum_{i=0}^{p-2} \delta (x)^i\delta (y)^{p-2-i}\delta (xy)$$
and again $\delta(x(y-x))=\delta (xy)$ such that this lies in 
$\lf (I)$. For $|x|=|y|$ even we have
$$q(x)-q(y)=q(x-y)-\delta (\sum_{i=1}^{p-1} \frac 1 i x^iy^{p-i})$$
so it suffices to see that $y-x$ divides the sum inside the
$\delta (-)$. The following equation in $\FF_p [x,y]$ 
shows that this is the case:
$$xy(y-x)\sum_{k=0}^{p-3}a_k x^ky^{p-3-k}=
\sum_{i=1}^{p-1} \frac 1 i x^i y^{p-i}\text{ where } 
a_k = \sum_{j=0}^k \frac 1 {j+1}.$$
The equation holds since by Euler's sum formula
$\sum_{n=1}^{p-1} n=0$ modulo $p$.
\end{proof}

\begin{definition}
Let $\kerd ,\imd ,\hd :\F \to \Alg$ denote the functors given by 
$$\kerd (R) = \ker (d), \quad \imd (R) = \image (d), \quad
\hd (R) = \kerd (R) / \imd (R)$$
where $d$ is the differential on $\drf (R)$.
\end{definition}

Recall from \cite{BO}, \cite{O1} and \cite{O2} that there are natural 
transformations of functors $\Phi : \til \to \hd$ and $Q: \lf \to \kerd$. 
It was shown that if $A\in \F$ is a free object, or its underlying  
algebra is polynomial, then $\Phi_A$ and $Q_A$ are isomorphisms.
We can now give a nice interpretation of the functor $\lf$, which
was originally defined by generators and complicated relations.

\begin{theorem}
\label{goodlf}
For any $R\in \F$ one has $\lf (R) \cong H_0(R; \kerd )$.
\end{theorem}

\begin{proof}
The induced map $Q_*:H_*(R;\lf ) \to H_*(R;\kerd )$ is an isomorphism since
$Q$ is an isomorphism on free objects. The 0th derived functor of 
$\lf$ was computed in Proposition \ref{0thapl}. 
\end{proof}

For any functor $E: \F \to \Alg$ we have that $H_i(A;F)=0$ for $i>0$ when
$A$ is a free object since we can use the trivial almost free 
resolution to compute the derived functors. 
For polynomial algebras we also have nice results. 

\begin{theorem}
\label{derivedwhenpoly}
Assume that the underlying algebra of $A \in \F$ is a polynomial
algebra. Then one has $H_i(A;\drf )=0$, $H_i(A;\til )=0$,
$H_i(A;\lf )=0$ and $H_i(A; \ell )=0$ for each $i>0$.
\end{theorem}

\begin{proof}
We first prove the statements for $\drf$ and $\til$.
Let $\Omega : \Alg \to \Alg$ denote the usual de Rham complex functor.
Pick an almost free resolution $P \in s\F$ of $A$. 
The forgetful functor $U:\F \to \Alg $ takes free objects to 
free objects. So we can apply $U$ to $P$ and get an almost free resolution
of $U(A)$ in $s\Alg$. Thus there is an isomorphism
$$H_i^\F (A;\Omega U) \cong H_i^\Alg (U(A);\Omega )$$
and the last group is trivial for $i>0$ since $U(A)$ is a free object in 
$\Alg$.

There is a linear map $\Omega U (A) \to \drf (A)$; 
$x_0 dx_1 \dots dx_n \mapsto x_0 dx_1 \dots dx_n$. The map is not 
multiplicative and it does not commute with the de Rham differential,
but it is an isomorphism of graded vector spaces.
Thus $H_i(A;\drf )$ is additively isomorphic to $H_i(A; \Omega U)$ which
is trivial for $i>0$. A similar isomorphism gives the result for 
the functor $\til$.

Next we consider the functor $\lf$.
The short exact sequence $0 \to \kerd \to \drf \to \imd \to 0$ gives
a long exact sequence of derived functors. By the above this 
sequence breaks up into the exact sequence
$$0 \to H_1(A;\imd ) \to H_0(A;\kerd ) \to H_0(A;\drf ) \to 
H_0(A;\imd ) \to 0$$
together with the isomorphisms $H_i(A;\kerd )\cong H_{i+1}(A;\imd )$
for $i\geq 1$.

There is also a short exact sequence $0 \to \imd \to \kerd \to \hd \to 0$
with corresponding long exact sequence of derived functors.
Since $\Phi$ is an isomorphism on free objects we have a
natural isomorphism $\Phi_*:H_*(-;\til ) \cong H_*(-;\hd )$.
By the above vanishing result for $H_*(A;\til )$ the long exact sequence 
breaks up into the short exact sequence 
$$0 \to H_0(A;\imd  ) \to H_0(A;\kerd ) \to H_0(A;\hd ) \to 0$$
and the isomorphisms $H_i(A;\imd )\cong H_i(A;\kerd )$ for $i\geq 1$.

Using Proposition \ref{0thapl} and Theorem \ref{goodlf} we can rewrite 
the exact sequences involving 0th derived functors as
$$0\to H_1(A;\imd )\to \lf (A) \to \drf (A) \to H_0(A;\imd ) \to 0$$
$$0\to H_0(A;\imd )\to \lf (A) \to \hd (A) \to 0.$$
Since $Q: \lf (A) \to \kerd (A)$ is an isomorphism we see that 
$H_1(A;\imd )=0$. By the isomorphisms 
$H_1(A;\imd )\cong H_1(A;\kerd )\cong H_2(A;\imd )\cong \dots$
we conclude that $H_i(A;\kerd )$ is trivial for $i>0$. 
But $H_*(-;\lf )$ is isomorphic to $H_*(-;\kerd )$ so we are done.

Finally we consider the functor $\ell$. By definition of $\lf$ there is
a short exact sequence $0 \to u\ell \to \ell \to \lf \to 0$. From 
the corresponding long exact sequence of derived functors we find that
$H_i(A;u\ell )\cong H_i(A;\ell )$ for $i>0$.
By Theorem \ref{filtration} from the appendix  
and Proposition 4.4 from \cite{O2} there is a short exact sequence
\begin{equation}
\label{sestilde}
0 \to u^{j+1}\ell (B) \to u^j\ell (B) \to u^j \otimes \til (B) \to 0
\quad , \quad j>0
\end{equation}
when $B$ is a free object in $\F$ or when the underlying algebra
of $B$ is a polynomial algebra.  
The corresponding long exact sequence of derived functors shows that
$H_i(A;u^j\ell ) \cong H_i(A;u^{j+1}\ell )$ so we have
$H_i(A;\ell ) \cong H_i(A; u^j \ell )$ for all $j\geq 0$.
But $(u^j\ell )^k=0$ for $k<2j$ so $H_i(A; u^j \ell)^k=0$ for $k<2j$
and the result follows. 
\end{proof}

\begin{proposition}
\label{0thell}
If the underlying algebra of an object $A\in \F$ is a polynomial 
algebra, then $H_0(A;\ell ) \cong \ell (A)$.
\end{proposition}

\begin{proof}
The short exact sequence $0 \to u\ell \to \ell \to \lf \to 0$
gives a short exact sequence of $0$th derived functors since
$H_1(A;\lf )=0$. Furthermore, there is a natural map
$H_0(-;F)\to F$ for any functor $F:\F \to \Alg$. 
So we have a commutative diagram with exact rows as follows:
$$
\begin{CD}
0 @>>> H_0(A;u\ell ) @>>> H_0(A;\ell ) @>>> H_0(A;\lf ) @>>> 0 \\
@. @VVV @VVV @VVV @. \\
0 @>>> u\ell (A) @>>> \ell (A) @>>> \lf (A) @>>> 0.
\end{CD}
$$
The right vertical map is an isomorphism so it suffices to 
show that the left vertical map is also an isomorphism.

Since $H_1(A;\til )=0$
the short exact sequence (\ref{sestilde}) gives a commutative diagram
as follows for $j>0$:
$$
\begin{CD}
0 @>>> H_0(A;u^{j+1}\ell ) @>>> H_0(A;u^j\ell ) 
@>>> H_0(A;u^j\otimes \til ) @>>> 0 \\
@. @VVV @VVV @VVV @. \\
0 @>>> u^{j+1} \ell (A) @>>> u^j\ell (A) @>>> u^j\otimes \til (A) @>>> 0.
\end{CD}
$$
where the right vertical map is an isomorphism.
Fix a degree $n$. For $j+1>n/2$ the map $H_0(A; u^{j+1}\ell )^n \to
(u^{j+1}\ell (A))^n$ is an isomorphism since both domain and target
space are zero. The result follows by induction. 
\end{proof}

\begin{corollary}
Let $X$ be a 1-connected space such that $H_*X$ is of finite type and
$H^*X$ is a polynomial algebra. Then then the spectral sequences of 
Theorem \ref{mainss} collapses at the $E_2$ terms. So there are 
isomorphisms
$$H_*(H^*(X) ;\drf )^* \cong H^*(\la X) \text{ and }
H_*(H^*(X) ;\ell )^* \cong H^*((\la X)_{h\TT }).$$
\end{corollary}

\section{The derived functors of an exterior algebra}

In the rest of this paper we take $p=2$.
Let $\ea =\Lambda (\sigma ) \in \F$ be an exterior algebra on one generator
of degree $|\sigma |\geq 2$. Note that $\lambda \sigma =0$ 
for dimensional reasons. We intend to compute the higher derived
functors of the various functors we have been considering for this
algebra. 

\begin{proposition}
\label{QuillenCor}
There are isomorphisms 
$$
H_*(\ea ;\drf )\cong \drf (\ea ) \otimes \Gamma [\omega ] \quad , \quad
H_*(\ea ;\til ) \cong \til (\ea ) \otimes \Gamma [\widetilde \omega].
$$
The inner degrees are $|\gamma_i(\omega ) |=i(2|\sigma |-1)$, 
$|\gamma_i(\widetilde \omega) |=i(4|\sigma |-1)$ and the grading
of the homology groups are given by
$$H_i(\ea ;\drf ) \cong \drf (\ea ) \otimes \gamma_i(\omega ) \quad , \quad
 H_i(\ea ;\til ) \cong \til (\ea ) \otimes \gamma_i(\widetilde \omega ).$$
\end{proposition} 

\begin{proof}
The algebra $\ea$ is the pushout of 
$\FF_2 \leftarrow \FF_2 [y] \to \FF_2 [x]$ where $y\mapsto x^2$.
Put $\lambda x=0$ and $\lambda y=0$. By Proposition \ref{QuillenSS}
we find
$$H_i(\ea ; E)\cong \Tor_i^{E(\FF_2 [y])}(\FF_2 , E(\FF_2 [x]))
\quad \text{for} \quad E=\drf , \til.$$
The result follows by standard computations.
\end{proof}

In order to compute derived functors of the other functors we need an
explicit simplicial resolution of $\ea$. By Theorem \ref{genres},
Theorem \ref{derivedwhenpoly}, Proposition \ref{0thapl} and
Proposition \ref{0thell} we may use an almost free resolution 
of $\ea$ in $s\Alg$ and equip it with $\lambda =0$.  

\begin{proposition}
\label{simplicial}
There is an almost free resolution $R_\bullet \in s\Alg $ of the
algebra $\ea$ with
$R_n= \FF_2 [x,y_1,y_2,\dots ,y_n]$ for $n\geq 0$.
The structure maps 
$d_i:R_n\to R_{n-1}$ and $s_i:R_n\to R_{n+1}$ 
are given by
\begin{eqnarray*}
s_i(x) &=& x \\
s_i(y_j) &=&
\begin{cases}
  y_j & ,j \leq i\\
  y_{j+1} &, j > i\\
\end{cases} \\
d_i(x) &=& x\\
d_i(y_j) &=&
\begin{cases}
  x^2 & , i=0, j=1\\
 y_{j-1}& , i <j, j>1\\ 
 y_{j} & , i \geq j, j<n \\
 0 & , i=n, j=n
\end{cases}
\end{eqnarray*}
The degrees of the generators are
$|x|=|\sigma|$ and $|y_i|=2|\sigma |$ for all $i$.
\end{proposition}

\begin{proof}
We first give a description of the simplicial set 
$\Delta_\bullet^1=\Hom_\scat (-,[1])$ suited for our purpose.

Define the elements $y_j\in \Delta^1_n$ 
for $n\geq 0$ and $0\leq j \leq n+1$ by
$y_j(i)=0$ if $i<j$ and $y_j(i)=1$ if $i\geq j$.
We have $\Delta_n^1 = \{ y_0, \dots , y_{n+1} \}$.
The structure maps are as follows:
$$
d_iy_j=
\begin{cases}
y_{j-1} & , i<j \\
y_j & , i\geq j
\end{cases}
\quad \text{and} \quad
s_iy_j=
\begin{cases}
y_{j+1} & , i<j \\
y_j & , i\geq j
\end{cases}
$$

Let $\FF_2 [-]$ denote the functor which takes a graded set
to the polynomial algebra generated by that set. 
Let $\FF_2 [\Delta^1 \simp , *]$ denote the pushout of
$\FF_2 \leftarrow \FF_2 [a] \to \FF_2 [\Delta^1 \simp ]$
where $\FF_2$ and $\FF_2[a]$ are constant simplicial algebras.
In degree $n$ the maps are as follows:
$a\mapsto 0 \in \FF_2$ and $a\mapsto y_{n+1}\in \FF_2 [\Delta^1\simp ]$.
Note that $\FF_2 [\Delta^1 \simp ] \simeq \FF_2 [*]$ by the simplicial
contraction of $\Delta^1 \simp $.
The spectral sequence \cite{Q} II \S 6 Theorem 6 gives that 
$\pi_i (\FF_2 [\Delta^1\simp ,*])\cong \FF_2$ for $i=0$ and $0$ otherwise. 

Define $R\simp$ as the pushout of 
$\FF_2 [x] \leftarrow \FF_2 [z] \to \FF_2 [\Delta^1\simp ,*]$ 
where in degree $n$ the maps are $z\mapsto x^2$ and $z\mapsto y_0$.
For this pushout Quillen's spectral sequence gives that
$$
\pi_i (R\simp )\cong \Tor_i^{\FF_2 [z]}(\FF_2 ,\FF_2 [x])= 
\begin{cases}
\ea & , i=0 \\
0 &, i>0
\end{cases}
$$
Thus $R\simp$ has the right homotopy groups. Further $R_n$ is as stated
and the structure maps are as stated. Note that $R\simp $ is almost
free. The degrees are correct since the structure maps must be 
degree preserving.
\end{proof}

\begin{lemma}
\label{identifygamma}
$H_n(\ea ;\Omega )$ has the following $\FF_2$-basis:
$$dy_1 \dots dy_n,\, xdy_1 \dots dy_n, \,
dxdy_1 \dots dy_n,\, xdxdy_1 \dots dy_n.$$
\end{lemma}

\begin{proof}
Using the formulas in Proposition \ref{simplicial} it is easy to check that 
the four given classes are in the kernel of $d_i$ for all $i$. 
To check linear independency,
we introduce two gradings of $\Omega $.

Firstly, the wedge grading on $\Omega (R_n)$ is defined as the number 
of wedge factors, ie. the number of $d$'s in a homogeneous element. 
Secondly the polynomial grading is defined as follows: Give $x$ grading 1,   
$y_j$ grading 2 and each $dx$ or $dy_j$ grading 0 and extend
multiplicatively. Note that the maps $d_i$ preserve both gradings.
We write $\Omega^{q,t} (R_n)$ for the elements in $\Omega (R_n)$ of wedge
degree $q$ and polynomial degree $t$. Thus there is a direct 
sum decomposition
$$H_n(\ea ;\Omega )= \bigoplus_{q,t\geq 0} H_n(\ea ;\Omega^{q,t} ).$$

The classes we consider sit in different bigradings, so we
only have to check that they individually do not represent the
trivial class. 

We have the following bases for $\Omega^{n,0} (R_{n+1})$, 
$\Omega^{n,0} (R_{n})$ and $\Omega^{n,0} (R_{n-1})$ respectively:

\begin{align*}
& \{ dy_1 \dots \widehat{dy_j} \dots dy_{n+1} \} \cup
\{dx dy_1 \dots \widehat{dy_j} \dots \widehat{dy_k} \dots dy_{n+1} \} , \\
& \{ dxdy_1 \dots \widehat{dy_j} \dots dy_n \} \cup
\{ dy_1 \dots dy_n \} , \\
&  \{ dxdy_1 \dots dy_{n-1} \}. 
\end{align*}

We use the normalized complex consisting of $\cap_{i>0} \ker (d_i)$ 
with differential $d_0$ to compute the homology. 
For this normalized complex we have the respective bases
$
\emptyset ,
\{ dy_1 \dots dy_n \} ,
\{ dx dy_1 \dots dy_{n-1} \} 
$.
Taking homology and using that $\Omega^{n,0} (R_{n-2})=0$ we see that  
the classes 
$dxdy_1\dots dy_{n-1}$ and
$dy_1\dots dy_{n}$ do not represent zero.

Similarly, $x dxdy_1\dots dy_{n-1}$ and
$x dy_1\dots dy_{n}$ do not represent zero. 
Keeping track of degrees and dimensions, the result follows.
\end{proof}

\begin{lemma}
\label{identifygammatilde}
$H_n(\ea ;\hd )$ has an $\FF_2$-basis as follows:
\begin{align*}
& y_1\dots y_{n}dy_1\dots dy_{n}, \, 
x^2y_1\dots y_{n} dy_1\dots dy_{n}, \\
& xdxy_1\dots y_{n}dy_1\dots dy_{n}, \,
x^3dxy_1\dots y_{n}dy_1\dots dy_{n}.
\end{align*}
\end{lemma}

\begin{proof}
This follows from Lemma \ref{identifygamma} and the Cartier isomorphism.
\end{proof}

\begin{corollary}
\label{split}
The long exact sequence associated to $\imd \to \kerd \to \hd$ 
splits into short exact sequences when evaluated at $\ea$:
$$0 \to H_*(\ea ;\imd ) \to H_*(\ea ; \kerd ) \to H_*(\ea ; \hd ) \to 0.$$
\end{corollary}

\begin{proof}
The generators from Lemma \ref{identifygammatilde} are
visibly in the image of the map 
$H_n(\ea ;\kerd ) \to H_n(\ea ;\hd )$.
\end{proof}

\begin{definition}
Let $\gamma_{i,j}\in \Omega  (R_j)$ denote the following elements:
\begin{align*}
& \gamma_{i,j}=y_1y_2\dots y_idy_1dy_2\dots dy_j 
\quad \text{for} \quad i,j>0 \\
& \gamma_{0,j}=dy_1dy_2\dots dy_j \quad \text{for} \quad j>0 \quad 
\text{and} \quad \gamma_{0,0} =1. 
\end{align*}
\end{definition}

These elements are actually in $\kerd (R_j)$ if $i\leq j$, and
in $\imd (R_j)$ in case $i <j$. They are even in
the normalized chain complex, so they define
classes in $H_*(\ea ;\kerd)$ respectively $H_*(\ea ;\imd)$.

\begin{theorem}
\label{determineB}
For $n\geq 0$ the homology group 
$H_{2n}(\ea ;\imd )$ has $\FF_2$-basis
$$\{ dx \gamma_{0,2n} \} \cup 
\{ \gamma_{2i,2n} \, ,\, x^2 \gamma_{2i,2n} \, , \,
xdx\gamma_{2i,2n} \, ,\, x^3dx\gamma_{2i,2n} | 0\leq i < n \} $$
and the homology group 
$H_{2n+1}(\ea ;\imd )$ 
has $\FF_2$-basis
\begin{align*}
& \{ \gamma_{0,2n+1} \, ,\, dx\gamma_{0,2n+1} \, ,\, xdx\gamma_{0,2n+1} \} 
\cup \\ 
& \{ \gamma_{2i+1,2n+1} \, , \, x^2 \gamma_{2i+1,2n+1} \, , \,
xdx\gamma_{2i+1,2n+1} \, ,\, x^3dx\gamma_{2i+1,2n+1} | 0\leq i < n \}.
\end{align*}
Similarly, $H_{2n}(\ea ;\kerd )$ has $\FF_2$-basis
$$\{ dx \gamma_{0,2n} \} \cup 
\{ \gamma_{2i,2n} \, ,\, x^2 \gamma_{2i,2n} \, , \,
xdx\gamma_{2i,2n} \, ,\, x^3dx\gamma_{2i,2n} | 0\leq i \leq n \} $$
and $H_{2n+1}(\ea ;\kerd )$ has $\FF_2$-basis
\begin{align*}
& \{ \gamma_{0,2n+1} \, ,\, dx\gamma_{0,2n+1} \, ,\, xdx\gamma_{0,2n+1} \} 
\cup \\ 
& \{ \gamma_{2i+1,2n+1} \, , \, x^2 \gamma_{2i+1,2n+1} \, , \,
xdx\gamma_{2i+1,2n+1} \, ,\, x^3dx\gamma_{2i+1,2n+1} | 0\leq i \leq n \}.
\end{align*}
\end{theorem}

\begin{proof}
Recall that $H_m(\ea ;\Omega )\cong \Omega (\ea ) \otimes \gamma_{0,m}$
and $H_m(\ea ;\hd ) \cong \til (\ea ) \otimes \gamma_{m,m}$.
From the  splitting
$H_m(\ea ;\kerd ) \cong H_m(\ea ;\imd ) \oplus 
H_m(\ea ;\hd )$
together with the computation of $H_m(\ea ;\hd )$ 
in Lemma \ref{identifygammatilde} ,
it follows that the statements about $H_m(\ea; \kerd)$ and
$H_m(\ea; \imd)$ are equivalent.

Recall the long exact sequence 
\begin{equation}
\label{lesext}
\begin{CD}
H_{m+1}(\ea ;\imd ) @>b>> H_m(\ea ;\kerd ) @>i_*>> 
H_m(\ea ;\Omega ) @>d_*>> H_m(\ea ; \imd ) 
\end{CD} 
\end{equation}
\emph{Claim:} The image of $d_*$ is a one dimensional
vector space, generated by a class, which in the normalized
chain complex is represented by  $dx \gamma_{0,m}$. 

Note that three out of four of the generators of 
$H_m(\ea;\Omega)$ are annihilated by $d_*$.  
So the image of $d_*$ is spanned by the single class,
represented in the normalized complex by
$dx\gamma_{0,m}\in \imd(\ea)$. 
This element actually represents a non zero homology class, since
the map $H_m(\ea ;\imd ) \to H_m(\ea ;\Omega )$
induced by inclusion, send it to the element represented by
$dxdy_1\dots dy_m$. But this is non-trivial according to 
Lemma \ref{identifygamma}, and the claim follows.

We will now prove the theorem by induction on $m$.
We first treat the case $m=0$.
Here we see that the map
$d_0:H_0(\ea ;\Omega )\to H_0(\ea ;\imd )$
is surjective by the long exact sequence (\ref{lesext}), 
and the statement follows directly from the claim we just proved.

Assume that the theorem holds for $m$. 
The long exact sequence (\ref{lesext}) gives us a short exact 
sequence
$$
0 \to \image (d_{m+1}) \to H_{m+1}(\ea;\imd) \overset b\to
\ker (i_m) \to 0
$$
So, in order to prove the theorem, we have to show that 
$b$ maps the classes given in the statement of the theorem
\emph{with $dx\gamma_{0,m+1}$ removed} to a basis for
the kernel of $i_m$.

By induction we have a basis for $H_{m}(\ea ;\kerd )$. 
From this we see that the kernel of the map 
$i_{2n+1}:H_{2n+1}(\ea ;\kerd ) \to H_{2n+1}(\ea ;\Omega )$
has basis
$$\{ \gamma_{2i+1,2n+1} \, ,\, x^2\gamma_{2i+1,2n+1} \, ,\,
xdx\gamma_{2i+1,2n+1} \, ,\, x^3dx\gamma_{2i+1,2n+1} |0\leq i \leq n \}$$
and that the kernel of the map
$i_{2n}:H_{2n}(\ea ;\kerd ) \to H_{2n}(\ea ;\Omega )$
has basis
$$\{ x^2 \gamma_{0,2n} \, , \,  x^3 dx\gamma_{0,2n}\} \cup 
\{ \gamma_{2i,2n} \, ,\, x^2 \gamma_{2i,2n} \, , \,
xdx\gamma_{2i,2n} \, ,\, x^3dx\gamma_{2i,2n} | 1\leq i \leq n \} .$$

It is convenient now to pass from the normalized to the
unnormalized complex. The normalized complex is
a subcomplex of the unnormalized one, and in the notation
we do not distinguish between a class in the subcomplex and
its image in the unnormalized one. 

Let us first consider the odd case, $m=2n+1$.
The element $\gamma_{2i,2n+2} \in \imd (R_{2n+2})$ is a cycle with
respect to the boundary $\partial = \sum d_k$. 
We compute its image under the map $b$. 

Let $\alpha_r \in \Omega (R_{2n+2})$ denote the following element:
$$
\alpha_r =y_1y_2\dots y_{2i}y_rdy_1dy_2\dots \widehat{dy_r} \dots dy_{2n+2}
\quad , \quad 2i+1\leq r \leq 2n+2
$$
where the hat means that the factor is left out. Put 
$\beta = \sum_{r=2i+2}^{2n+2} \alpha_r$.
We have $d\alpha_r=\gamma_{2i,2n+2}$ for each $r$, so that
also $d\beta =\gamma_{2i,2n+2}$.
This means that $b(\gamma_{2i,2n+2})$ is represented by
$$
\partial \beta =
\sum_{r=2i+2}^{2n+2}{(y_{r-1}+y_r)\gamma_{2i,2n+1}}+y_{2n+1}\gamma_{2i,2n+1}
=\gamma_{2i+1,2n+1}.
$$
Since $b$ is linear with respect to multiplication by 
$x^2$, $xdx$, $x^3dx$ this gives the desired result for
$H_{2n+2}(\ea ;\imd )$. 

In the even case $m=2n$, a similar argument shows that 
$b(\gamma_{2i+1,2n+1})$ is represented by $\gamma_{2i+2,2n}$. 

Checking with the lists of classes above, we see that we 
are left to prove that $b(\gamma_{0,2n+1})=x^2\gamma_{0,2n}$.
The argument is very similar. Let 
$$\alpha_r=y_rdy_1dy_2\dots \widehat{dy_r}\dots dy_{2n+1} \text{ and } 
\beta =\alpha_1 +\alpha_2 +\dots + \alpha_{2n+1}.$$
Then $d\alpha_r =\gamma_{0,2n+1}$ so also $d\beta_r =\gamma_{0,2n+1}$.
Thus $b(\gamma_{0,2n+1})$ is represented by 
$\partial \beta = x^2 \gamma_{0,2n}$.
\end{proof}

\begin{proposition}
There are short exact sequences for $i\geq 0$, $t\geq 1$ as follows:
\begin{align*}
& \begin{CD}
0 @>>> H_i(\ea ;u\ell ) @>>> H_i(\ea ;\ell ) @>>> H_i(\ea ; \lf ) @>>> 0 \, , 
\end{CD}
\\
& \begin{CD}
0 @>>> H_i(\ea ;u^{t+1}\ell ) @>>> H_i(\ea ;u^t\ell ) @>>> 
H_i(\ea ; u^t\til ) @>>> 0. 
\end{CD}
\end{align*}
\end{proposition}

\begin{proof}
The first short exact sequence follows if we can prove that the
connecting homomorphism $b:H_{i+1}(\ea ;\lf )\to H_{i}(\ea ;u\ell )$ is
trivial. So consider the diagram
$$
\begin{CD}
0 @>>> u\ell (R_{i+1}) @>>> \ell (R_{i+1}) @>>> \lf (R_{i+1}) @>>> 0 \\
@. @VVV @VVV @VVV @. \\
0 @>>> u\ell (R_{i}) @>>> \ell (R_{i}) @>>> \lf (R_{i}) @>>> 0. 
\end{CD}
$$

The element $q(y_1) \dots q(y_j) \delta (y_{j+1}) \dots \delta (y_{i+1})
\in \ell (R_{i+1})$ maps to the element $\gamma_{j,i+1} \in \lf (R_{i+1})$.
By the relations $\delta (a)^2 =\delta (\lambda a)$, 
$q(a)^2=\phi (\lambda a)+\delta (a^2\lambda a)$ and 
$\delta (a)q(b)= \delta (a\lambda b)+ \delta (ab) \delta (b)$ we see that
this element maps down to zero in $\ell (R_i)$. So the connecting 
homomorphism $b$ is trivial.

The proof for the second short exact sequence is similar. 
\end{proof}

\begin{corollary}
\label{lderived}
For each $i\geq 0$ there are isomorphisms of $\FF_2$-vector spaces 
\begin{align*}
H_i(\ea ;\ell ) & \cong H_i(\ea ;\lf )\oplus
\big( \bigoplus_{t\geq 1} u^t\otimes H_i(\ea ;\til ) \big) \\
& \cong H_i(\ea ;\imd )\oplus
\big( \FF_2 [u]\otimes H_i(\ea ;\til ) \big).
\end{align*}
\end{corollary}

Let $F:\F \to \Alg$ be a functor. We define the total
degree of a class in $H_i(R;F)^n$ to be $n-i$. For $F=\ell, \drf$ this
corresponds through the spectral sequences of Theorem \ref{mainss} to the 
grading of cohomology groups.
We write $\PS_F(t)$ for the Poincar\'e series corresponding to
the total degree of $H_*(\ea ;F)^*$.

\begin{theorem}
\label{Poincare}
Let $\sdeg $ denote the degree of the class $\sigma\in \ea$.  
Then we have the following Poincar\'e series: 
\begin{align*}
\PS_\drf (t)&=(1+t^{\sdeg })(1-t^{\sdeg-1})^{-1}, \\ 
\PS_\til (t) &=(1+t^{2\sdeg})(1-t^{2\sdeg-1})^{-1}, \\
\PS_\imd (t) &= t^{\sdeg -1}(1+t^{\sdeg -1}-t^{2\sdeg -1})
(1-t^{2\sdeg -1})^{-1}(1-t^{2\sdeg -2})^{-1}, \\
\PS_\ell(t) &=(1+t^{\sdeg -1}-t^{\sdeg +1}+t^{2\sdeg -1})
(1-t^2)^{-1}(1-t^{2\sdeg -2})^{-1}. 
\end{align*}
\end{theorem}

\begin{proof} 
The first two formulas follow from
Proposition \ref{QuillenCor}: The total degree of 
$\gamma_i (\omega )$ is $|\gamma_i (\omega )|-i=(2\sdeg -2)i$
and $\drf (\ea )= \ea \otimes \Lambda (d\sigma )$ so
$$
\PS_\drf (t)=(1+t^{\sdeg })(1+t^{\sdeg -1})(1-t^{2\sdeg -2})^{-1}
=(1+t^\sdeg )(1-t^{\sdeg -1})^{-1}.
$$
A similar argument gives $\PS_\til (t)$.

To determine $\PS_\imd (t)$ we must count 
the classes given in Theorem \ref{determineB}, according to the 
total degree.  

We divide these classes into three groups. The first group are
those of the form $dx \gamma_{0,m}$. The Poincar\'e
series of the subspace generated by those classes is 
$t^{\sdeg -1}/(1-t^{2\sdeg-2})$. The second group are those 
of the type $\gamma_{0,2n+1}$ or $xdx\gamma_{0,2n+1}$.
These have Poincar\'e series 
$t^{2\sdeg -2}(1+t^{2\sdeg -1})/(1-t^{4\sdeg -4})$. 

The third group is the remaining classes. They span a free 
$\til (\ea )$-module with basis 
$X=\{ \gamma_{2i,2n} ,\gamma_{2i+1,2n+1} | 0\leq i <n, 0\leq n \}.$ 
We introduce the following operation on the set $X$:
$T(\gamma_{i,n})=\gamma_{i+1,n+1}$. This operation has total
degree $4\sdeg-2$. All generators are obtained by applying $T$
a non-negative number of times starting from one of the elements of 
$Y=\{\gamma_{0,2n} |n\geq 1 \}$.

The set $Y$ has Poincar\'e series $t^{4\sdeg -4}/(1-t^{4\sdeg-4})$.
So, the set $X$ has Poincar\'e series
$t^{4\sdeg -4}(1-t^{4\sdeg-4})^{-1}(1-t^{4\sdeg -2})^{-1}$. We 
multiply this by $(1+t^{2\sdeg })(1+t^{2\sdeg -1})$, make a small 
reduction, and obtain
the Poincar\'e series for the third group of classes:
$$t^{4\sdeg -4}(1+t^{2\sdeg })(1-t^{2\sdeg -1})^{-1} (1-t^{4\sdeg -4})^{-1}.$$
To get the Poincar\'e series of $\imd$, we add the three series 
obtained so far. 

$$
\PS_\imd (t)= \frac {t^{\sdeg -1}} {1-t^{2\sdeg -2}}+
\frac {t^{2\sdeg -2}(1+t^{2\sdeg -1})} {1-t^{4\sdeg -4}}+
\frac {t^{4\sdeg -4}(1+t^{2\sdeg })} {(1-t^{2\sdeg -1})(1-t^{4\sdeg -4})}
$$
The stated formula for $\PS_\imd (t)$ follows after some reductions.

Finally, Corollary \ref{lderived} gives that
$\PS_\ell (t)=\PS_\imd (t)+(1-t^2)^{-1}\PS_\til (t)$ which leads to
the stated formula after some reductions.
\end{proof}

\section{The spectral sequences for spheres}

Let $X$ be a pointed space, $Y=\Sigma X$ it's reduced suspension.
We have established a spectral sequence converging to
$H^*(E\TT \times_{\TT } \Lambda Y;\FF_2)$ in general. But in this
special case, we are fortunate to have a direct calculation of
the homology $H_*(E\TT_+\wedge_{\TT } \Lambda Y;\FF_2)$. If the homology
of $X$ is of finite type, the (finite) dimensions of these homology and
cohomology groups agree. So if we for the particular space $X$ can
check that the Poincar\'e series of the homology 
of $H^*(E\TT \times_{\TT } \Lambda Y;\FF_2)$ as computed in \cite{CC} 
agrees with the Poincar\'e series of the $E_2$ term of our spectral
sequence, we know that our spectral sequence collapses. 
(See also \cite{BM} for an easier proof of the results in \cite{CC}).

Let us now consider the special case of spheres $X=S^{s-1}$ and
$Y=S^s$. The purpose of this section is to show that in this case the
two Poincar\'e series actually agree, forcing the spectral sequence to
collapse.

\begin{theorem}
\label{PSbyCC}
The Poincar\'e series of 
$H^*(E\TT \times_{\TT } \Lambda S^s;\FF_2 )$ is
$$
(1+t^{\sdeg-1}-t^{\sdeg+1}+t^{2\sdeg-1})(1-t^{2})^{-1}(1-t^{2\sdeg-2})^{-1}.
$$
\end{theorem}

\begin{proof}
We compute a sequence of related Poincar\'e series.
First, let $A=\tilde H^*(S^{s-1})$ considered as a graded
vector space. This has Poincar\'e series $t^{s-1}$.

For each $m\geq 1$, the cyclic group $C_m$ acts on $A^{\otimes m}$.
This is a 1-dimensional vector space over $\FF_2$. 
We now consider the homology groups
$$
H_*(C_m;A^{\otimes m})
$$
We first look at homological dimension 0. 
$H_0(C_m,A^{\otimes m})\cong A^{\otimes m}$, so it has Poincar\'e series
$t^{m(s-1)}$.
In higher homological degrees, there are two cases. 
If $m$ is odd, the groups all vanish, and we get a trivial Poincar\'e
series.
If $m$ is even, and $i \geq 1$,
$H_i(C_m,A^{\otimes m})\cong A^{\otimes m}$ Since this single group has
homological degree $i$, its Poincar\'e series is
$t^{i+m(s-1)}$.

Now, recall from \cite{CC}, proposition 9.3 that 
\[
\tilde H_*(ES^1_+\wedge_{S^1} \Lambda S^s) \cong
\oplus_{m\geq 1} H_*(C_m;A^{\otimes m})
\]
(Actually, we are correcting a misprint in \cite{CC} here. 
The homology groups on the right hand side of the formula 
should not be reduced).

The Poincar\'e series of the right hand side contains the sum of
the contribution of the homology in dimension zero. The
Poincar\'e series of this part is 
$\sum_{m\geq 1} {t^{m(s-1)}}=t^{s-1}(1-t^{s-1})^{-1}$.
It also contains the sums of the contributions of the
reduced group homologies. Since this is trivial if $m$ is even,
we can as well put $m=2n$, and the Poincar\'e series of the
reduced part is 
$$
\sum_{i\geq 1}\sum_{n\geq 1} t^{i+2n(\sdeg-1)}=
t^{2(\sdeg-1)+1}(1-t)^{-1}(1-t^{2(\sdeg-1)})^{-1}
$$
Summing, we get that the Poincar\'e series for
$\tilde H_*(ES^1_+\wedge_{S^1} \Lambda(\Sigma X))$ is
$$
(t^{\sdeg-1}-t^\sdeg+t^{2\sdeg-2})(1-t)^{-1}(1-t^{2\sdeg-2})^{-1}
$$
Finally, we note that there is a short exact sequence of
homology groups:
$$
0 \to \tilde H_*(B\TT ) \to \tilde H_*(E\TT \times_{\TT } \Lambda S^s)
\to \tilde H_*(E\TT_+ \wedge_{\TT } \Lambda S^s) \to 0.
$$
This shows that the Poincar\'e series of 
$H_*(E\TT \times_{\TT } \Lambda S^s)$ is 
$$
(t^{\sdeg-1}-t^\sdeg+t^{2\sdeg-2})(1-t)^{-1}(1-t^{2\sdeg-2})^{-1}
+(1-t^2)^{-1}
$$
Bringing on common denominator and adding proves the theorem.
\end{proof}

\begin{proposition}
\label{PSbyKY}
The Poincar\'e series of $H^*(\Lambda S^s;\FF_2 )$ is
$(1+t^s)(1-t^{s-1})^{-1}$ when $s\geq 2$.
\end{proposition}

\begin{proof}
The mod 2 cohomology ring of $\Lambda S^s$ is a special case of 
Theorem 2.2 of \cite{KY} except for the case $s=2$. It is however
shown (Remark 2.6) that the Eilenberg-Moore spectral sequence also 
collapses when $s=2$ so we can compute the Poincar\'e series from 
the $E_2$-term. It has the following form 
(see the proof of Theorem 2.2):
$$E_2^{*,*} \cong \Lambda (x) \otimes \Lambda (\overline x) \otimes
\Gamma [\omega ]$$
where the respective bidegrees of $x$, $\overline x$ and 
$\gamma_i (\omega )$ are
$(0,s)$, $(-1,s)$ and $(-2i,2is)$ 
such that the respective total degrees becomes
$s$, $s-1$ and $2i(s-1)$.
Thus the Poincar\'e series is
$$(1+t^s)(1+t^{s-1})(1-t^{2(s-1)})^{-1}$$
and the result follows by a small reduction.
\end{proof}

\begin{theorem}
If we let $X=S^s$ with $s\geq 2$ and use $\FF_2$-coefficients, then
the spectral sequences of Theorem \ref{mainss} collapses.  
Thus there are isomorphisms of graded $\FF_2$-vector spaces:
$$H_*(H^*(S^s);\drf )^* \cong H^*(\Lambda S^s) \text{ and }
H_*(H^*(S^s);\ell )^* \cong H^*((\Lambda S^s)_{h\TT }).$$ 
\end{theorem}

\begin{proof}
By Theorem \ref{Poincare}, Theorem \ref{PSbyCC} 
and Proposition \ref{PSbyKY} the Poincar\'e series
of the $E_2$-terms agree with the Poincar\'e series of the 
targets. So the spectral sequences collapses.
\end{proof}

\section{Appendix: On a filtration of the functor $\ell $}
\label{filtrationap}

In this appendix we identify the graded object associated with 
the filtration
$$\ell (A) \supseteq u\ell (A) \supseteq u^2\ell (A)\supseteq \dots$$
in the case where $p=2$ and $A$ is a polynomial algebra.

Recall that the functors $\lf, \til : \F \to \Alg$ are defined by
$\lf (A)= \ell (A)/ (u)$ and $\til (A)= \lf (A)/I_\delta (A)$ where 
$I_\delta (A)$ is the ideal $(\delta (x)|x\in A)\subseteq \lf (A)$.

We want to define a map $\ell(A) \to \til(A)[t]$ such that the elements
$\phi (x)$, $q(x)$ and $u$ in the domain are send to the elements 
$\phi (x)$, $q(x)$ and $t^2$ in the target. Unfortunately, this cannot
be done by a ring map. But if we pay the penalty of changing the
multiplicative structure of the target, we can almost get such a map.

\begin{definition}
$\Ot (A)$ is the free graded commutative algebra on
generators $\phi (x), q(x)$ for $x\in A$ and $t$, of degrees 
$|\phi (x)|=2|x|$, $|q(x)|=2|x|-1$ and $|t|=1$, modulo
the relations
\begin{align*}
& q(x+y)=q(x)+q(y), & & \phi (x+y)= \phi (x)+\phi (y), \\ 
& q(xy)=\phi (x)q(y)+\phi (y)q(x), & & \phi (xy)=\phi (x)\phi (y), \\ 
& q(x)^2=\phi (\lambda x)+tq(\lambda x).  
\end{align*}
\end{definition}

Clearly, $\Ot (A)/(t)\cong \til (A)$. The ring $\Ot (A)$ is just a 
twisted version of the polynomial ring over $\til (A)$ in $t$ in the
following case:

\begin{theorem}
\label{gradedot}
Assume that the underlying algebra of $A$ is a polynomial algebra.  
Then the graded ring $Gr_*(\Ot (A))$ corresponding to the filtration
of $\Ot (A)$ by powers of $t$ equals $\til (A)[t]$.
\end{theorem}

\begin{proof}
As an intermediate step, let us consider the ring $R(A)$ which is
defined exactly like $\Ot (A)$ except that we do not include the
last relation $q(x)^2=\phi (\lambda x)+tq(\lambda x)$.

If $A$ is a polynomial algebra
on generators $\{ x_i|i\in I \}$, then $R(A)$ is a polynomial ring on 
generators $\phi (x_i)$ and $q(x_i)$. To obtain $\Ot (A)$ from $R(A)$,
we have to add the relations 
$q(p)^2=\phi (\lambda p)+tq(\lambda p)$, where $p$
is any polynomial in the generators $x_i$. Actually,
it is sufficient to do this for the generators themselves, as this
relation for $p_1p_2$ follows from the relations for $p_1$ and $p_2$.
Because, assume those are satisfied, then we calculate
\begin{align*}
q(p_1p_2)^2 &= \phi (p_1)^2q(p_2)^2+\phi (p_2)^2q(p_1)^2 \\
&= \phi (p_1)^2 (\phi (\lambda p_2)+tq(\lambda p_2))+
\phi (p_2)^2(\phi (\lambda p_1)+tq(\lambda p_1)) \\
&= \phi (\lambda (p_1p_2))+tq(\lambda (p_1p_2)).
\end{align*}

Thus we can write $\Ot (A)$ as an algebra:
$$
\FF_2 [t,\phi (x_i), q(x_i)|i\in I ]/
\{ q(x_i)^2 = \phi (\lambda x_i)+tq(\lambda x_i) \}
$$
From this it is clear, that $\Ot (A)$ is a free 
$\FF_2 [t, \phi (x_i)|i\in I]$-module, with
generators 
$$\{ q(x_{i_1}) \dots q(x_{i_n})|i_r\neq i_s \text{ for } r\neq s, 
\, n\geq 0\}.$$
(The empty product means 1 here.) It follows that
$Gr_*(\Ot (A))$ is a free module over $Gr_* (\FF_2 [t,\phi (x_i)|i\in I] )$ 
with the same generators.

So, to finish the proof, we only have to determine the multiplicative
structure of $\Ot (A)$. The multiplicative relations are given 
by the relations. In the graded ring they are $q(x_i)^2=\phi (\lambda x_i)$.
So, we have a presentation of the graded ring as
$$
\FF_2 [t,\phi (x_i),q(x_i)|i \in I]/ \{ q(x_i)^2=\phi(\lambda x_i) \} .
$$
But this is exactly  $\til (A)[t]$.
\end{proof}

\begin{theorem}
\label{filtration}
Let $A$ be an object in $\F$ and $i\geq 1$ an integer. Multiplication 
with $u^i$ defines a natural surjective $\FF_2$-linear map
$$u^i: \til (A) \to \frac {u^i\ell (A)} {u^{i+1}\ell (A)}.$$
If the underlying algebra of $A$ is a polynomial algebra,
then this map is an isomorphism and  
$$Gr_*(\ell (A))\cong \lf (A)\oplus \bigoplus_{j\geq 1} u^j\otimes \til (A).$$
\end{theorem}

\begin{proof}
Multiplication with $u^i$ gives a surjective map $\ell (A) \to u^i\ell (A)$
and $u^i(I_\delta (A))=0$, $u^i(u\ell (A))=u^{i+1}\ell (A)$ so the map 
factors through $\til (A)$.

We define a natural ring map $v: \ell (A) \to \Ot (A)$ by the formulas
$$v(\phi (x))=\phi (x)+tq(x),\quad v(q(x))=q(x),\quad v(u)=t^2,\quad 
v(\delta (x))=0.$$
To see that $v$ is well defined, we have to check that
the relations in the definition of $\ell$ goes to 0. This is trivial
for all relations except three which is verified as follows:
\begin{align*}
v(\phi (xy))&= \phi (xy)+tq(xy)=(\phi (x)+tq(x))(\phi (y)+tq(y))+t^2q(x)q(y)\\
            &= v(\phi(x)\phi(y)+uq(x)q(y)),\\
v(q(xy))&= q(xy)=\phi (x)q(y)+q(x)\phi (y)\\
&= (\phi (x)+tq(x))q(y)+(\phi (y)+tq(y))q(x)\\
&= v(\phi (x)q(y)+q(x)\phi(y)),\\  
v(q(x)^2) &= q(x)^2=v(\phi (\lambda x)+\delta (x^2\lambda x)). 
\end{align*}

By the map $v$ we get a commutative diagram as follows:
$$
\begin{CD}
\til (A) @>>> \Ot (A)/t^2 \Ot (A)\\
@V{u^i}VV                                @V{t^{2i}}VV     \\
u^i\ell (A)/u^{i+1}\ell (A) @>>> t^{2i} \Ot (A)/t^{2i+2} \Ot (A)\\
\end{CD}
$$
When the underlying algebra of $A$ is a polynomial algebra, then 
Theorem \ref{gradedot} gives that the top and the right vertical maps 
are injective. So in this case the left vertical map is also injective. 
\end{proof}

\end{document}